\documentclass[11pt]{article}
\usepackage{amsmath,amssymb}

-\headheight\advance\topmargin by -\headsep

\def\N{\bf \mbox{I\hspace{-.15em}N}}
\def\Z{\bf \mbox{Z\hspace{-.40em}Z}}
\def\R{{\bf \mbox{I\hspace{-.20em}R}}}

\def\bkR{{\bf \mbox{I\hspace{-.20em}R}}}
\def\C{C^{\infty}(M, {\bf \mbox{I\hspace{-.20em}R}})}
\def\lcf{\lbrack\! \lbrack}
\def\rcf{\rbrack\! \rbrack}

\newtheorem{definition}{Definition}[section]
\newtheorem{lemma}[definition]{Lemma}
\newtheorem{proposition}[definition]{Proposition}
\newtheorem{theorem}[definition]{Theorem}
\newtheorem{remark}[definition]{Remark}
\newtheorem{corol}[definition]{Corollary}

\newtheorem{examples}[definition]{Examples}
\newtheorem{example}[definition]{Example}
\newenvironment{proof}{\noindent{\bf Proof.}}{\hfill $\square$}
\newenvironment{proof-of}{\noindent{\bf Proof of}}{\hfill $\square$}

\begin{document}

\title{Twisted Jacobi manifolds, twisted Dirac-Jacobi structures and quasi-Jacobi bialgebroids}
\author{J. M. Nunes da Costa \\Departamento de Matem\'atica \\ Universidade de Coimbra\\Apartado 3008
\\3001-454 Coimbra - Portugal\\ {\small E-mail: jmcosta@mat.uc.pt}
\and F. Petalidou\\Faculty of Sciences and Technology
\\University of Peloponnese \\22100 Tripoli - Greece\\{\small E-mail: petalido@uop.gr}}

\date{}
\maketitle

\begin{abstract}

We study twisted Jacobi manifolds, a concept that we had introduced
in a previous Note. Twisted Jacobi manifolds can be characterized
using twisted Dirac-Jacobi, which are sub-bundles of Courant-Jacobi
algebroids. We show that each twisted Jacobi manifold has an
associated Lie algebroid with a $1$-cocycle. We introduce the notion
of quasi-Jacobi bialgebroid and we prove that each twisted Jacobi
manifold has a quasi-Jacobi bialgebroid canonically associated.
Moreover, the double of a quasi-Jacobi bialgebroid is a
Courant-Jacobi algebroid. Several examples of twisted Jacobi
manifolds and twisted Dirac-Jacobi structures are presented.

\end{abstract}

\vspace{3mm} \noindent {\bf{Keywords : }}{Twisted Jacobi manifold,
twisted Dirac-Jacobi structure, Jacobi bialgebroid, Courant-Jacobi
algebroid,
 quasi-Jacobi bialgebroid.}

\vspace{3mm} \noindent {\bf A.M.S. classification (2000):} 53D10,
53D17, 17Bxx.

\section{Introduction}
Jacobi manifolds were introduced by Lichnerowicz \cite{lch} and
Kirillov \cite{kr} as smooth manifolds endowed with a bivector field
$\Lambda$ and a vector field $E$ satisfying some compatibility
conditions. When the vector field $E$ identically vanishes, the
Jacobi manifold is just a Poisson manifold. So, Poisson manifolds
are particular cases of Jacobi manifolds. But there are other
examples of Jacobi structures on manifolds which are not Poisson,
such as contact structures  and local conformally symplectic
structures.

The notion of twisted Poisson manifold (or Poisson manifold with a
$3$-form background) was introduced by \v{S}evera and Weinstein
\cite{sw}, motivated by the works of Klim\v{c}ik and Strobl
\cite{kl} on topological field theory and Park \cite{p} on string
theory. Since Jacobi structures on manifolds generalize Poisson
structures, the introduction of the concept of a twisted Jacobi
manifold seems very natural. This task was achieved in the Note
\cite{jf} where, besides we have introduced that notion, we briefly
presented some of its properties.

Dirac structures on manifolds were introduced by Courant and
Weinstein \cite{cw} and developed in detail by Courant \cite{c}.
Dirac structures include presymplectic forms, Poisson structures and
foliations. The first approach to extend the theory of Dirac
structures to Jacobi manifolds was done by Wade \cite{wd}, who
introduced the $\mathcal{E}^1(M)$-Dirac structures as a natural
extension of Dirac bundles in the sense of Courant \cite{c}. These
$\mathcal{E}^1(M)$-Dirac structures, which we call Dirac-Jacobi
structures, include Jacobi manifolds and are sub-bundles of the
vector bundle $(TM \times \R) \oplus (T^*M \times \R)$ over $M$,
satisfying a certain integrability condition. However, the vector
bundle $(TM \times \R) \oplus (T^*M \times \R)$ is not a Courant
algebroid. This fact motivated a more general treatment, proposed in
\cite{gm2,jj}. The concept of Courant-Jacobi algebroid was
introduced, independently, in \cite{gm2} and \cite{jj}, and the main
example of this structure is the double of a Jacobi bialgebroid
\cite{im1, gm1}. A Dirac structure for a Courant-Jacobi algebroid
$E$ is defined as a sub-bundle of the vector bundle $E$ over $M$
satisfying an integrability condition. Dirac-Jacobi bundles arise
then as a particular case of these structures.

As we have already mentioned, twisted Poisson manifolds were
introduced by \v{S}evera and Weinstein \cite{sw} who studied them in
the framework of Courant algebroids and Dirac structures. For the
case of twisted Jacobi  manifolds, we use Dirac-Jacobi structures.
More precisely, we use twisted Dirac-Jacobi structures, which are
sub-bundles of the Courant-Jacobi algebroid $(TM \times \R) \oplus
(T^*M \times \R)$ equipped with a ``twisted bracket" on its space of
sections. These  Dirac-Jacobi bundles enable us to characterize
twisted Jacobi structures on manifolds.

On the other hand, Roytenberg \cite{ro} developed a theory of
quasi-Lie bialgebroids and used it to study twisted Poisson
manifolds \cite{ro1}. Namely, with each twisted Poisson structure on
a manifold $M$, a quasi-Lie bialgebroid structure on $(TM, T^*M)$
can be associated. When we try to investigate what happens in the
Jacobi framework, we realize that things are different. First of all
because, in opposition to the Poisson case, one cannot, in general,
define a Lie algebroid structure on the cotangent bundle $T^*M$ of a
Jacobi manifold  $(M,\Lambda,E)$. Usually, only the vector bundle
$T^*M \times \R$ over $M$ admits such a structure \cite{krb}.
Furthermore, with each Jacobi manifold, there exists an associated
Jacobi bialgebroid \cite{im1, gm1}, while in the case of a Poisson
manifold it admits an associated Lie bialgebroid. Motivated by these
facts, we introduce the concept of a quasi-Jacobi bialgebroid, which
is the one that fits in our theory. We prove that each twisted
Jacobi manifold has an associated quasi-Jacobi bialgebroid and that
the double of a quasi-Jacobi bialgebroid is a Courant-Jacobi
algebroid.

The paper is divided into eight sections. In section 2 we recall
some facts on Jacobi manifolds and their relation with Lie algebroid
theory. In section 3 we study the main properties of a twisted
Jacobi manifold, we present some examples and we show that if $M$ is
equipped with a twisted Jacobi structure, then there exists a
twisted exact homogeneous Poisson structure on $M \times \R$.
Section 4 is devoted to twisted Dirac-Jacobi structures and we
characterize twisted Jacobi manifolds using these structures.
Several examples of twisted Dirac-Jacobi bundles are presented,
including graphs of sections of $\bigwedge^2(T^*M \times \R)$ and
twisted locally conformal presymplectic structures. We also relate
twisted Dirac-Jacobi bundles and Dirac bundles in the sense of
Courant. In Section 5 we see how gauge transformations act on
twisted Dirac-Jacobi structures. In section 6 we construct a Lie
algebroid with a $1$-cocycle associated with each twisted Jacobi
manifold. The notion of quasi-Jacobi bialgebroid is introduced in
section 7 and we prove that its double is a Courant-Jacobi
algebroid. In section 8 we show that each twisted Jacobi manifolds
admits an associated quasi-Jacobi bialgebroid.

\vspace{3mm} \noindent {\bf Notation :} In this paper, $M$ is a
$C^{\infty}$-differentiable manifold of finite dimension. We denote
by $TM$ and $T^{\ast}M$, respectively, the tangent and cotangent
bundles over $M$ and by $C^{\infty}(M,\R)$ the space of all real
$C^{\infty}$-differentiable functions on $M$. For the Schouten
bracket and the interior product of a form with a multivector field,
we use the convention of sign indicated by Koszul \cite{kz}, (see
also \cite{mrl}).

\section{Jacobi manifolds}
A {\em Jacobi manifold} is a differentiable manifold $M$ equipped
with a bivector field $\Lambda$ and a vector field $E$ such that
\begin{equation} \label{1.1}
[\Lambda,\Lambda]=-2E\wedge \Lambda
\hspace{5mm}\mathrm{and}\hspace{5mm}[E,\Lambda] = 0,
\end{equation}
where $[\cdot,\cdot]$ denotes the Schouten bracket \cite{kz}. In
this case, $(\Lambda,E)$ defines a bracket on $C^{\infty}(M,\R)$
which is called the {\em Jacobi bracket} and is given, for all
$f,g\in C^{\infty}(M, \R)$, by
\begin{equation} \label{1.2}
\{f,g\} = \Lambda (df,d g) +  f(E.g) - g (E.f).
\end{equation}
The Jacobi bracket endows $C^{\infty}(M,\R)$ with a local Lie
algebra structure in the sense of Kirillov \cite{kr}. Reciprocally,
a local Lie algebra structure on $C^{\infty}(M, \R)$ induces on $M$
a Jacobi structure.

When the vector field $E$ identically vanishes on $M$, the Jacobi
structure reduces to a Poisson structure on the manifold. However,
there are other examples of Jacobi manifolds either than Poisson
manifolds, such as contact and locally conformal symplectic
manifolds, \cite{lch}.

There are some well-known results concerning Jacobi structures on
manifolds that we briefly recall.

Let $(M,\Lambda,E)$ be a Jacobi manifold. Then, the pair
$(\Lambda,E)$ defines the homomorphism of $\C$-modules
$(\Lambda,E)^{\#} : \Gamma(T^{*}M \times \R) \rightarrow \Gamma(TM
\times \R)$ given, for any section $(\alpha,f)$ of $T^{*}M \times
\R$, by
\begin{equation} \label{1.3}
(\Lambda,E)^{\#}(\alpha,f) = (\Lambda^{\#}(\alpha) + fE, - i_E
\alpha),
\end{equation}
and, with each $f\in \C$, we can associate the vector field
$X_f=\Lambda^{\#}(df)+fE$, called the \emph{hamiltonian vector
field} of $f$. We have that
$$
X_f=\pi ((\Lambda,E)^{\#}(df,f)),
$$
where $\pi : TM\times \R \to TM$ denotes the projection over the
first factor. Moreover, for all $f,g \in \C$,
\begin{equation} \label{1.4}
[X_f,X_g]=X_{\{f,g\}}.
\end{equation}
Also, the vector bundle $T^*M\times \R$ over $M$ endowed with the
anchor map $\pi \circ (\Lambda,E)^\# : T^*M\times \R \to TM$ and the
Lie algebra bracket $\{\cdot,\cdot\}$ on the space of its sections,
given, for all $(\alpha,f),(\beta,g) \in \Gamma(T^*M\times \R)$, by
\begin{equation}\label{1.5}
\{(\alpha,f),(\beta,g)\} = (\gamma,r),
\end{equation}
where
$$
\gamma = {\cal L }_{\Lambda^{\#}(\alpha)}\beta - {\cal L
}_{\Lambda^{\#}(\beta)}\alpha - d(\Lambda(\alpha, \beta)) + f {\cal
L }_E\beta - g {\cal L }_E\alpha - i_E(\alpha \wedge \beta),
$$
$$
r = - \Lambda(\alpha, \beta) + \Lambda(\alpha,dg) -
\Lambda(\beta,df) +  fE(g) - gE(f),
$$
is a Lie algebroid over $M$ \cite{krb}. The associated exterior
derivative $d_*$ on $\Gamma(\bigwedge(TM\times \R))= \oplus_{k \in
\Z} \Gamma(\bigwedge^k(TM \times \R))$ is given \cite{lmp}, for all
$(P,Q) \in \Gamma(\bigwedge^k (TM\times \R))\cong \Gamma
(\bigwedge^k (TM))\oplus \Gamma (\bigwedge^{k-1} (TM))$, by
\begin{equation} \label{1.6}
d_{*} (P,Q)= ([\Lambda,P]+ k E \wedge P+  \Lambda \wedge Q,
-[\Lambda,Q]+(1-k)E \wedge Q+[E,P]).
\end{equation}

\vspace{1mm}

It is well known that, given a Lie algebroid $(A,[\cdot,\cdot], a)$
over a differentiable manifold $M$ with a $1$-cocycle $\phi \in
\Gamma(A^*)$ in the Lie algebroid cohomology complex with trivial
coefficients \cite{mck}, we can modify the usual representation of
the Lie algebra $(\Gamma(A), [\cdot,\cdot])$ on $\C$ by defining a
new representation $a^{\phi}:\Gamma (A) \times \C \rightarrow \C$ as
\begin{equation}\label{1.12}
a^{\phi}(X,f)= a(X)f+ (i_X \phi) f, \quad \quad  \forall \, (X,f)
\in \Gamma(A) \times \C.
\end{equation}
Therefore, we obtain a new cohomology operator $d^\phi$ on
$\Gamma(\bigwedge A^*)= \oplus_{k \in \Z} \Gamma(\bigwedge^k A^* )$
given by
\begin{equation} \label{1.10}
d^{\phi}(\beta)=d \beta+ \phi \wedge \beta, \quad \quad  \forall \,
\beta \in \Gamma({\bigwedge}^k A^*),
\end{equation}
where $d$ is the cohomology operator defined by $([\cdot,\cdot],a)$
on $\Gamma(\bigwedge A^*)$, and a new Lie derivative operator of
forms with respect to $X\in \Gamma(A)$, ${\mathcal L}^{\phi}_{X}=
d^\phi \circ i_X + i_X \circ d^\phi$, that can be expressed in terms
of the usual Lie derivative $\mathcal{L}_X = d\circ i_X + i_X \circ
d$, as
\begin{equation}\label{1.11}
{\mathcal L}^{\phi}_{X}(\beta)={\mathcal L}_X \beta+ (i_X \phi)
\beta, \quad \quad \forall \, \beta \in \Gamma({\bigwedge}^k A^*).
\end{equation}
Using $\phi$, it is also possible to modify the Schouten bracket
$[\cdot,\cdot]$ on the graded algebra $\Gamma(\bigwedge A)=\oplus_{k
\in \Z}\Gamma(\bigwedge^k A)$ to the $\phi$-Schouten bracket
$[\cdot,\cdot]^{\phi}$ on $\Gamma(\bigwedge A)$ defined, for any $P
\in \Gamma({\bigwedge}^p A)$ and $Q \in \Gamma({\bigwedge}^q A)$, by
\begin{equation} \label{defor}
[P,Q]^{\phi}=[P,Q]+(p-1)P \wedge (i_{\phi}Q)+(-1)^p (q-1)(i_{\phi}P)
\wedge Q,
\end{equation}
where $i_{\phi}Q$ and $i_{\phi}P$ can be interpreted as the usual
contraction of a multivector field with a $1$-form. A differential
calculus using $a^{\phi}$, $d^{\phi}$, ${\mathcal L}^{\phi}$ and
$[\cdot,\cdot]^{\phi}$ can be developed. The formulae obtained are
similar, but adapted, to the case of a Lie algebroid \cite{im1},
\cite{gm1}.

A pair $(A,\phi)$ formed  by a Lie algebroid $A$ and a $1$-cocycle
$\phi$ of $A$, is called a {\em Jacobi algebroid} in the terminology
of \cite{gm1}.

\vspace{3mm}

A trivial example of a Jacobi algebroid over $M$ is the vector
bundle $TM\times \R \to M$ equipped with the bracket
\begin{equation} \label{1.7}
[(X,f),(Y,g)] = ([X,Y],X. g - Y. f), \quad \quad \forall \,
(X,f),(Y,g) \in \Gamma(TM\times \R),
\end{equation}
the vector bundle map $\pi : TM \times \R \to TM$, that is the
projection over the first factor, and the section $(0,1)$ of
$T^*M\times\R$. The associated exterior derivative on
$\Gamma(\bigwedge (T^*M\times \R))$ is the operator ${\rm d}=(d,
-d)$ and $(0,1)$ is a $1$-cocycle in the cohomology complex with
trivial coefficients of $(TM\times \R, [\cdot,\cdot],\pi,
\mathrm{d})$. In the sequel, we will denote by $\mathrm{d}^{(0,1)}$
the differential operator on $\Gamma(\bigwedge (T^*M\times \R))$
modified by $(0,1)$, as in (\ref{1.10}).

\vspace{3mm}

The notion of \emph{generalized Lie bialgebroid} and the equivalent
one of \emph{Jacobi bialgebroid} were introduced, respectively, by
D. Iglesias and J.C. Marrero in \cite{im1} and by J. Grabowski and
G. Marmo in \cite{gm1} in such a way that a Jacobi manifold has a
Jacobi bialgebroid canonically associated and conversely. A Jacobi
bialgebroid over $M$ is a pair $(A,A^*)$ of Lie algebroids over $M$,
in duality, with differentials $d$ and $d_*$, respectively, endowed
with a $1$-cocycle $\phi\in \Gamma(A^*)$ of $(A,d)$ and a
$1$-cocycle $W\in \Gamma(A)$ of $(A^*,d_*)$, such that, for every $P
\in \Gamma (\bigwedge^{p}A)$ and $Q \in \Gamma (\bigwedge A)$, the
following condition holds :
$$
d_{*}^{W} [P,Q]^{\phi}= [d_{*}^{W} P,Q]^{\phi} + (-1)^{p+1}[P,
d_{*}^{W}Q]^{\phi}.
$$

\vspace{3mm}

The pair formed by the Jacobi algebroid $(TM\times \R,
[\cdot,\cdot],\pi,(0,1))$, presented above, together with the Lie
algebroid $(T^*M\times \R, \{\cdot,\cdot\}, \pi \circ
(\Lambda,E)^{\#})$ and the $1$-cocycle $(-E,0) \in \Gamma(TM\times
\R)$ on it, is a Jacobi bialgebroid over the Jacobi manifold
$(M,\Lambda,E)$, \cite{im1}.

\vspace{3mm}

Finally, let us recall \cite{im1} that a section $(\Lambda,E)$ of
$\bigwedge^2(TM \times \R)$ defines a Jacobi structure on the
manifold $M$ if and only if
\begin{equation} \label{1.9}
[(\Lambda, E), (\Lambda, E)]^{(0,1)}=(0,0).
\end{equation}

\section{Twisted Jacobi manifolds}
In \cite{jf} we introduced the concept of twisted Jacobi manifold
and we presented some of its properties. Now, in this section, we
will review and complete the results announced in \cite{jf}.

\vspace{3mm}

We start by recalling that, given a bivector field $\Lambda$ on a
differentiable manifold $M$, the associated vector bundle map
$\Lambda^{\#} : T^*M \to TM$ induces a homomorphism of $\C$-modules
$\Lambda^{\#} : \Gamma(T^*M) \to \Gamma(TM)$,
$$
\langle \beta, \Lambda^{\#}(\alpha)\rangle = \Lambda(\alpha,\beta),
\quad \quad \forall \, \alpha, \beta \in \Gamma(T^*M),
$$
that can be extended to a homomorphism, also denoted by
$\Lambda^{\#}$, from $\Gamma ({\bigwedge}^k( T^*M))$ onto $\Gamma
({\bigwedge}^k(TM ))$, $k\in \N$, as follows:
\begin{equation}\label{ext-hom}
\Lambda^{\#}(f) = f  \quad \mathrm{and} \quad (\Lambda^{\#}
\eta)(\alpha_1, \ldots, \alpha_k)=(-1)^k \eta( \Lambda^\#
(\alpha_1), \ldots, \Lambda^\# (\alpha_k)),
\end{equation}
for all $f\in \C$, $\eta \in \Gamma ({\bigwedge}^k( T^*M))$ and
$\alpha_1,\ldots,\alpha_k \in \Gamma(T^*M)$. Analogously, with each
section $(\Lambda,E)$ of $\bigwedge^2(TM\times \R)$, we can
associate a homomorphism of $\C$-modules
$$
(\Lambda,E)^\#:\Gamma ({\bigwedge}^k( T^*M \times \R)) \to \Gamma
({\bigwedge}^k(TM \times \R)), \quad k\in \N,
$$
by setting, for all $f\in \C$, $(\eta,\xi) \in \Gamma (\bigwedge^k(
T^*M \times \R))$ and $(\alpha_1,f_1), \ldots,(\alpha_k,f_k)\in
\Gamma(T^*M \times \R)$,
$$
(\Lambda,E)^\#(f)=f
$$
and
\begin{eqnarray} \label{2.10}
\lefteqn{(\Lambda,E)^\#(\eta,\xi)((\alpha_1,f_1), \ldots,
(\alpha_k,f_k)) }  \nonumber \\
& & =(-1)^k(\eta,\xi)((\Lambda,E)^\# (\alpha_1,f_1), \cdots ,
(\Lambda,E)^\# (\alpha_k,f_k)).
\end{eqnarray}
We remark that for $k=1$, we recover (\ref{1.3}).

\vspace{3mm}

Let us introduce some notation, following \cite{sw}. Let $\Lambda$
be a bivector field on $M$ and $\varphi$ a $3$-form on $M$. We
denote by $(\Lambda^{\#} \otimes 1)(\varphi)$ the section of
$(\bigwedge ^2 TM) \otimes T^*M$ that acts on multivector fields by
contraction with the factor in $T^*M$. For any $f\in \C$, $X \in
\Gamma(TM)$ and $\alpha, \beta \in \Gamma(T^*M)$,
\begin{equation} \label{2.3}
( \Lambda^{\#} \otimes 1)(\varphi)(f)=0 \quad \mathrm{and} \quad
(\Lambda^{\#} \otimes 1)(\varphi)(\alpha,\beta)(X)= -
\varphi(\Lambda^{\#}(\alpha), \Lambda^{\#}(\beta), X).
\end{equation}
Similarly, if $\omega$ is a $2$-form on $M$, then, for any $X \in
\Gamma(TM)$ and $\alpha \in \Gamma(T^*M)$,
$$
(\Lambda^{\#} \otimes
1)(\omega)(\alpha)(X)=\omega(\Lambda^{\#}(\alpha), X).
$$

\vspace{1mm}

In what follows, we consider the Jacobi algebroid $(TM\times
\R,[\cdot,\cdot],\pi,(0,1))$ and we are mainly interested in the
vector bundle map defined by (\ref{2.10}) for $k=3$.

\begin{proposition} \label{p2.1}
Let $(\Lambda,E)$ be a section of $\bigwedge^2(TM\times \R)$ and
$(\varphi,\omega)$ a section of $\bigwedge^3(T^*M\times \R)$. Then,
$$
[(\Lambda, E), (\Lambda, E)]^{(0,1)}= 2
(\Lambda,E)^{\#}(\varphi,\omega)
$$
if and only if
\begin{equation}\label{2.5}
[\Lambda,\Lambda] +2E \wedge \Lambda  =  2 \Lambda^{\#}(\varphi)+ 2(
\Lambda^{\#} \omega)\wedge E
\end{equation}
and
\begin{equation} \label{2.6}
[E,\Lambda]=  ( \Lambda^{\#} \otimes 1)(\varphi)(E)- (( \Lambda^{\#}
\otimes 1)(\omega)(E)) \wedge E.
\end{equation}
\end{proposition}
\begin{proof}
Let $(\alpha, f), (\beta,g), (\gamma,h)$ be three arbitrary sections
of $T^*M\times \R$. We have,
\begin{eqnarray}  \label{2.8}
\lefteqn{[(\Lambda, E), (\Lambda, E)]^{(0,1)}((\alpha, f),
(\beta,g), (\gamma,h))    \nonumber} \\
& = &([\Lambda, \Lambda]+2 E \wedge \Lambda, 2 [E,
\Lambda])((\alpha, f), (\beta,g), (\gamma,h)) \nonumber \\
& =& ([\Lambda,\Lambda]+ 2E \wedge \Lambda)(\alpha, \beta, \gamma)
+2f [E, \Lambda](\beta, \gamma)-2g [E, \Lambda](\alpha, \gamma) \nonumber \\
& & + 2h [E, \Lambda](\alpha, \beta).
\end{eqnarray}
On the other hand,
\begin{eqnarray} \label{2.9}
\lefteqn{2 (\Lambda,E)^{\#}(\varphi,\omega)((\alpha, f),
(\beta,g), (\gamma,h)) \nonumber }\\
& = & 2( \Lambda^{\#}\varphi)(\alpha, \beta, \gamma)+ 2((
\Lambda^{\#} \omega)\wedge E)(\alpha, \beta,
\gamma) \nonumber \\
& & - 2 \left(\varphi(\Lambda^{\#}(\beta), \Lambda^{\#}(\gamma),
fE)-\varphi(\Lambda^{\#}(\alpha), \Lambda^{\#}(\gamma),gE)+
\varphi(\Lambda^{\#}(\alpha), \Lambda^{\#}(\beta),hE) \right)
\nonumber  \\
& & - 2 ( (i_E \alpha)[\omega(\Lambda^{\#}(\gamma), gE)-
\omega(\Lambda^{\#}(\beta), hE)]-(i_E
\beta)[\omega(\Lambda^{\#}(\gamma), fE)
\nonumber \\
& & -\omega(\Lambda^{\#}(\alpha), hE)] +(i_E
\gamma)[\omega(\Lambda^{\#}(\beta), fE)-
\omega(\Lambda^{\#}(\alpha), gE)] ) \nonumber \\
& = & 2\left((\Lambda^{\#} \varphi)+(
(\Lambda^{\#} \omega)\wedge E), ( \Lambda^{\#} \otimes 1)(\varphi)(E)  \right. \nonumber \\
& & \left. - (( \Lambda^{\#} \otimes 1)(\omega)(E) \wedge E) \right)
((\alpha, f), (\beta,g), (\gamma,h)).
\end{eqnarray}
Comparing the terms on trivector fields and bivector fields of
(\ref{2.8}) and (\ref{2.9}), we obtain, respectively, the formul\ae
\, (\ref{2.5}) and (\ref{2.6}).
\end{proof}

\vspace{.3cm}

The sections of $\bigwedge^3(T^*M\times \R)$ that are closed with
respect to the differential operator $\mathrm{d}^{(0,1)}$ will have
a special role hereafter. We will call them $\mathrm{d}^{(0,1)}${\em
-closed}.

\begin{lemma}\label{lem-closed}
A section $(\varphi,\omega)$ of $\bigwedge^3(T^*M\times \R)$ is
$\mathrm{d}^{(0,1)}$-closed, i.e.
$\mathrm{d}^{(0,1)}(\varphi,\omega)=(0,0)$, if and only if $\varphi
= d\omega$.
\end{lemma}
Thus, we shall denote any $\mathrm{d}^{(0,1)}$-closed section
$(\varphi,\omega)$ of $\bigwedge^3(T^*M\times \R)$ by
$(d\omega,\omega)$, with $\omega$ a $2$-form on $M$.

\begin{definition}
A {\em twisted Jacobi structure} on a differentiable manifold $M$ is
defined by choosing a bivector field $\Lambda$, a vector field $E$
and a $2$-form $\omega$ on $M$ such that
\begin{equation}  \label{2.1}
[(\Lambda, E), (\Lambda, E)]^{(0,1)}= 2
(\Lambda,E)^{\#}(d\omega,\omega).
\end{equation}
A manifold equipped with such a structure is called a {\em twisted
Jacobi manifold} or a $ \omega${\em -Jacobi manifold} and it is
denoted by the triple $(M,(\Lambda,E),\omega)$.
\end{definition}

Hence, according to Proposition \ref{p2.1}, we may define a twisted
Jacobi manifold as a manifold $M$ equipped with a section
$(\Lambda,E)$ of $\bigwedge^2(TM \times \R)$ and a $2$-form $\omega$
on $M$ satisfying conditions (\ref{2.5}) and (\ref{2.6}), for
$\varphi = d\omega$.

\begin{examples}\label{exemples-tJ}
\end{examples}
\vspace{-2mm}

\noindent \emph{1. Jacobi manifolds:} Any Jacobi manifold
$(M,\Lambda,E)$ endowed with a $2$-form $\omega$ satisfying
$(\Lambda,E)^\#(d\omega,\omega)=(0,0)$ can be viewed as a twisted
Jacobi manifold.

\vspace{1mm}

\noindent \emph{2. Twisted locally conformal symplectic manifolds:}
A \emph{twisted locally conformal symplectic manifold} is a
$2n$-dimensional differentiable manifold $M$ equipped with a
non-degenerate $2$-form $\Theta$, a closed $1$-form $\vartheta$,
called the \emph{Lee $1$-form}, and a $2$-form $\omega$ such that
$$
d(\Theta + \omega) + \vartheta \wedge (\Theta + \omega) =0.
$$
Let $E$ be the unique vector field and $\Lambda$ the unique bivector
field on $M$ which are defined by
\begin{equation}\label{eq-tlcs}
i(E)\Theta = -\vartheta  \hspace{5mm} \mathrm{and} \hspace{5mm}
i(\Lambda^\#(\alpha))\Theta = -\alpha, \hspace{5mm} \mathrm{for}\,\,
\mathrm{all}\,\, \alpha \in \Gamma(T^*M).
\end{equation}
If we also denote by $\Lambda^\#$ the extension (\ref{ext-hom}) of
the isomorphism $\Lambda^\# : \Gamma(T^*M) \to \Gamma(TM)$ given by
(\ref{eq-tlcs}), we obtain
$$
E = \Lambda^\#(\vartheta) \hspace{3mm} \mathrm{and} \hspace{3mm}
\Lambda = \Lambda^\#(\Theta).
$$
By a simple, but very long computation, we prove that the pair
$((\Lambda,E),\omega)$ satisfies the relations (\ref{2.5}) and
(\ref{2.6}), for $\varphi=d\omega$. Whence, $((\Lambda,E),\omega)$
endows $M$ with a twisted Jacobi structure.

\vspace{1mm}

\noindent \emph{3. A trivial example in local coordinates:} Let
$(x_0,x_1,x_2,x_3,x_4)$ be a system of local coordinates in $\R^5$.
Let us consider a bivector field $\Lambda$, a vector field $E$ and a
2-form $\omega$ on $\R^5$ given, in these coordinates, by
$$
\displaystyle{\Lambda=\frac{\partial}{\partial x_1} \wedge
\frac{\partial}{\partial x_3} + \frac{\partial}{\partial x_2} \wedge
\frac{\partial}{\partial x_4} +x_4 \frac{\partial}{\partial x_0}
\wedge \frac{\partial}{\partial x_4}}, \,\,\,\,\,
\,\displaystyle{E=\frac{\partial}{\partial x_0}}, \,\,\,\,\,\,
\omega= dx_1 \wedge dx_3.
$$
A simple computation gives
$$
[\Lambda, \Lambda]+ 2 E \wedge \Lambda=
2\displaystyle{\frac{\partial}{\partial x_1} \wedge
\frac{\partial}{\partial x_3} \wedge \frac{\partial}{\partial x_0}}
\quad \quad \mathrm{and} \quad \quad [E, \Lambda]=0.
$$
Since
$$
\Lambda^\#(\omega)= \displaystyle{\frac{\partial}{\partial x_1}
\wedge \frac{\partial}{\partial x_3}} \quad \quad \mathrm{and} \quad
\quad (\Lambda^\# \otimes 1)(\omega)(E) =0,
$$
we have
$$
[\Lambda, \Lambda]+ 2 E \wedge \Lambda= 2\Lambda^\#(\omega) \wedge E
\quad \quad \mathrm{and} \quad \quad [E, \Lambda]=-(\Lambda^\#
\otimes 1)(\omega)(E) \wedge E.
$$
According to Proposition \ref{p2.1}, with $\varphi=d \omega=0$,
$((\Lambda,E),\omega)$ defines a twisted Jacobi structure on the
manifold $\R^5$.

\vspace{3mm}

Given a twisted Jacobi structure $((\Lambda,E),\omega)$ on $M$,
$(\Lambda,E)$ defines on $\C$ an internal composition law
$\{\cdot,\cdot\}$ just as in the case of Jacobi structure: For all
$f,g \in \C$,
\begin{equation}\label{br-M}
\{f,g\} = \Lambda (df,d g) +  fE(g) - g E(f).
\end{equation}
Since (\ref{1.9}) does not hold, this bracket fails the Jacobi
identity and is no more a Lie bracket.

\begin{proposition}  \label{p3.4}
Let $(M,(\Lambda,E),\omega)$ be a twisted Jacobi manifold. Then, for
all $f,g,h \in  \C$,
$$
\{f,\{g,h\}\}+ c.p.= -(d
\omega,\omega)((\Lambda,E)^{\#}(df,f),(\Lambda,E)^{\#}(dg,g),
(\Lambda,E)^{\#}(dh,h)),
$$
where c.p. denotes sum after circular permutation.
\end{proposition}
\begin{proof}
The result follows directly from (\ref{2.10}) for $k=3$ and
(\ref{2.1}), taking into account that, for any $f,g,h \in \C$,
$$
\frac{1}{2}[(\Lambda, E), (\Lambda,
E)]^{(0,1)}((df,f),(dg,g),(dh,h))=\{f,\{g,h\}\}+ p.c.
$$
\end{proof}

Let us now examine some relations between twisted Jacobi manifolds
and twisted Poisson manifolds.

\vspace{3mm}

We recall that a \emph{twisted Poisson manifold} \cite{sw} is a
differentiable manifold $M$ endowed with a bivector field $\Lambda$
and a closed $3$-form $\varphi$ on $M$ such that $[\Lambda,
\Lambda]= 2 \Lambda^\# (\varphi)$. When $\varphi$ is exact, i.e.
$\varphi = d\omega$ with $\omega \in \Gamma(\bigwedge^2T^*M)$, we
say that $(M,\Lambda,\varphi)$ is a {\em twisted exact} Poisson
manifold. A twisted Jacobi manifold $(M,(\Lambda,E),\omega)$, with
$E=0$, defines a twisted exact Poisson structure on $M$, since
$$
[(\Lambda,0), (\Lambda,0)]^{(0,1)}= 2 (\Lambda, 0)^\# (d \omega,
\omega)\,\,  \Rightarrow \,\, [\Lambda,\Lambda]=2 \Lambda ^\#(d
\omega).
$$

Furthermore, it is well known that there exists a close relationship
which links homogeneous Poisson manifolds with Jacobi manifolds
\cite{lch}. Namely, to each Jacobi manifold $(M,\Lambda,E)$,
 we can associate a homogeneous Poisson
manifold $(\tilde{M}, \tilde{\Lambda},
\displaystyle{\frac{\partial}{\partial t}})$, called the {\em
Poissonization} of $(M,\Lambda,E)$, with $\tilde{M} = M \times \R$
and $\tilde{\Lambda}= e^{-t}(\Lambda+
\displaystyle{\frac{\partial}{\partial t}} \wedge E)$,  $t$ being
the canonical coordinate on $\R$. For the twisted exact Poisson
structures, we introduce the following definition.

\begin{definition}  \label{def.tehp}
A {\em homogeneous twisted exact Poisson structure} on a manifold
$M$ is defined by a triple $(\Lambda, Z, \omega)$, where $\Lambda$
is a bivector field on $M$, $Z$ is a vector field on $M$ and
$\omega$ is a $2$-form on $M$, such that
$$
[\Lambda, \Lambda]= 2\Lambda^\# (d \omega), \quad [Z, \Lambda]= -
\Lambda,\quad L_Z\omega= \omega.
$$
\end{definition}

\begin{proposition}
Let $(M,(\Lambda,E), \omega)$ be a twisted Jacobi manifold. We set
$\tilde{M}= M\times \R$ and we consider on $\tilde{M}$ the tensor
fields $\tilde{\Lambda}= e^{-t}(\Lambda+
\displaystyle{\frac{\partial}{\partial t}} \wedge E)$ and
$\tilde{\omega}= \,e^t \omega$, $t$ being the canonical coordinate
on the factor $\R$. Then, the triple
$(\tilde{\Lambda},\displaystyle{\frac{\partial}{\partial
t}},\tilde{\omega})$ defines an homogeneous twisted exact Poisson
structure on $\tilde{M}$.
\end{proposition}
\begin{proof}
We have $[\displaystyle{\frac{\partial}{\partial
t}},\tilde{\Lambda}]=-\tilde{\Lambda}$ and $L_{\partial /
\partial t}\tilde{\omega}= \tilde{\omega}$. So, according to Definition
\ref{def.tehp}, it remains to prove that $[\tilde{\Lambda},
\tilde{\Lambda}]= 2\tilde{\Lambda}^\#(d \tilde{\omega})$. From the
definition of $\tilde{\Lambda}$, we compute
$$
[\tilde{\Lambda}, \tilde{\Lambda}]= e^{-2t}([\Lambda, \Lambda]+ 2 E
\wedge \Lambda) + 2e^{-2t}( \frac{\partial}{\partial t} \wedge [E,
\Lambda])
$$
and, since $(M,(\Lambda,E), \omega)$ is a twisted Jacobi manifold,
from (\ref{2.5}) and (\ref{2.6}), we can write
\begin{eqnarray} \label{3.30}
[\tilde{\Lambda}, \tilde{\Lambda}]& = & 2e^{-2t} \left( \Lambda^\#
(d \omega)  + \Lambda^\# (\omega) \wedge E \right. \nonumber  \\
& & \left. + \frac{\partial}{\partial t} \wedge \left((\Lambda^\#
\otimes 1) (d \omega) (E) - ((\Lambda^\# \otimes 1) ( \omega) (E))
\wedge E\right) \right).
\end{eqnarray}
On the other hand,
\begin{equation} \label{3.31}
\tilde{\Lambda}^\#(d \tilde{\omega})=  \, e^t \tilde{\Lambda}^\#(d
\omega + d t \wedge \omega).
\end{equation}
But,
\begin{equation} \label{3.32}
\tilde{\Lambda}^\#(d \omega)= e^{-3t} \left( \Lambda^\#(d
\omega)+\displaystyle{\frac{\partial}{\partial t}} \wedge
((\Lambda^\# \otimes 1)(d \omega)(E)) \right)
\end{equation}
and
\begin{equation} \label{3.33}
\tilde{\Lambda}^\#(d t \wedge \omega)=e^{-3t} \left(\Lambda^\#
(\omega)  - \displaystyle{\frac{\partial}{\partial t}} \wedge
((\Lambda^\# \otimes 1) ( \omega) (E)) \right) \wedge E.
\end{equation}
From equations (\ref{3.30})-(\ref{3.33}) we obtain
$[\tilde{\Lambda}, \tilde{\Lambda}]=2 \tilde{\Lambda}^\#(d
\tilde{\omega})$.
\end{proof}

\section{Twisted Dirac-Jacobi structures}
The notions of Courant-Jacobi algebroid and the equivalent one of
generalized Courant algebroid were introduced in \cite{gm2} and
\cite{jj}, respectively, as a generalization of the definition of
Courant algebroid \cite{lwx, ro}.

\begin{definition} {\rm(\cite{jj})} \label{def2}
A {\em generalized Courant algebroid} or a {\em Courant-Jacobi
algebroid} on a differentiable manifold $M$ is a vector bundle $E$
over $M$ equipped with a nondegenerate symmetric bilinear form
$(\cdot,\cdot)$ on the bundle, a skew-symmetric bracket
$[\cdot,\cdot]$ on $\Gamma(E)$, a bundle map $\rho^\theta: E \to
TM\times \R$ and a section $\theta$ of $E^*$ such that, for any
$e_1, e_2 \in \Gamma (E)$, the condition $\langle \theta
,[e_1,e_2]\rangle = \rho (e_1)\langle \theta, e_2 \rangle
-\rho(e_2)\langle \theta, e_1 \rangle $ holds, $\rho$ being the
bundle map from $E$ onto $TM$ induced by $\rho^\theta$, satisfying,
for all $e,e_1,e_2,e_3 \in \Gamma(E)$ the following properties:
\begin{itemize}
\item[i)] $[[e_1,e_2],e_3] + [[e_2,e_3],e_1]+[[e_3,e_1],e_2] =
{\cal D}^\theta T(e_1,e_2,e_3)$, \\where $T(e_1, e_2,
e_3)=\frac{1}{3} ([e_1, e_2], e_3) + c.p.$ and ${\cal D}^{\theta}:
C^{\infty}(M,\R) \to \Gamma(E)$ is the first-order differential
operator given by $ ({\cal D}^{\theta} f,e)= \frac{1}{2} \rho
^{\theta}(e)f$; \item[ii)] $\rho^\theta
([e_1,e_2])=[\rho^\theta(e_1),\rho^\theta(e_2)]$,\\
where the bracket on the right-hand side is the Lie bracket
(\ref{1.7}) on $\Gamma(TM \times \bkR)$; \item[iii)] $\rho^\theta
(e)(e_1,e_2) =  ([e,e_1]+ {\cal D}^{\theta} (e,e_1),e_2) +( e_1,
[e,e_2]+{\cal D}^{\theta} (e,e_2))$; \item[iv)] for any $f,g \in
C^\infty (M, \R)$, $ ({\cal D}^\theta f, {\cal D}^\theta g)=0$.
\end{itemize}
\end{definition}

A {\em Dirac structure} for the generalized Courant algebroid
$(E,\theta)$ is a sub-bundle $L$ of $E$ which is closed under the
bracket $[\cdot,\cdot]$ and is maximally isotropic with respect to
the symmetric bilinear form $(\cdot,\cdot)$. In this case $(L,
\rho|_L,[\cdot,\cdot]|_L)$ is a Lie algebroid over $M$.

\vspace{3mm}

An important example of a Courant-Jacobi algebroid is the double $A
\oplus A^*$ of a Jacobi bialgebroid $((A,\phi),(A^*,W))$ over $M$
\cite{gm2,jj}. The bracket on the space $\Gamma(A\oplus A^*)$ of its
sections is given, for all $e_1= X_1+\alpha_1, e_2= X_2+\alpha_2 \in
\Gamma(A\oplus A^*)$, by
\begin{eqnarray} \label{4.3}
\lcf X_1+\alpha_1,  X_2+\alpha_2 \rcf &=& \left ( [X_1,X_2]^{\phi}+
{\mathcal L}^{W}_{*{\alpha_1}}X_2 - {\mathcal
L}^{W}_{*{\alpha_2}}X_1 - d^{W}_{*}(e_1,e_2)_- \right ) \nonumber \\
& + & \left ([\alpha_1, \alpha_2]_{*}^{W} + {\mathcal
L}^{\phi}_{X_1}\alpha_2 - {\mathcal L}^{\phi}_{X_2}\alpha_1 +
d^{\phi}(e_1,e_2)_- \right),
\end{eqnarray}
where $(e_1,e_2)_- = \frac{1}{2}(i_{X_2}\alpha_1 -i_{X_1}\alpha_2)$.
Moreover, $\theta=\phi+W$, $\rho$ is the sum of the anchor maps of
$A$ and $A^*$, the symmetric bilinear form on $A\oplus A^*$ is the
canonical one, i.e. $(e_1,e_2)=(e_1,e_2)_+ =
\frac{1}{2}(i_{X_2}\alpha_1 +i_{X_1}\alpha_2)$, ${\mathcal D}=(d_*+
d)_{|C^{\infty}(M,\bkR)}$ and ${\mathcal D}^\theta = (d_* ^W+
d^\phi)_{|C^{\infty}(M,\bkR)}$.

For the case of the Jacobi bialgebroid $((TM\times \R, (0,1)),(T^*M
\times \R, (0,0)))$, where $T^*M \times \R$ is equipped with the
null Lie algebroid structure, the Courant-Jacobi structure defined
on its double ${\mathcal E}^1(M)= (TM\times \R) \oplus (T^*M \times
\R)$ corresponds to the following bracket on the space
$\Gamma({\mathcal E}^1(M))$, defined in \cite{wd} as a direct
generalization of the Courant bracket on $\Gamma(TM \oplus T^*M)$
\cite{c}, as follows: for all $e_1= (X_1,f_1)+(\alpha_1,g_1),\, e_2=
(X_2,f_2)+(\alpha_2,g_2) \in \Gamma({\mathcal E}^1(M))$,
\begin{eqnarray*}
\lcf e_1, e_2 \rcf  & = & \lcf (X_1,f_1)+(\alpha_1,g_1),
(X_2,f_2)+(\alpha_2,g_2) \rcf   \\
& = & [(X_1, f_1), (X_2,f_2)]^{(0,1)}  \\
&  &   +\left( {\mathcal L}^{(0,1)}_{ (X_1, f_1)} (\alpha_2,g_2) -
{\mathcal L}^{(0,1)}_{(X_2, f_2)} (\alpha_1,g_1)+
\mathrm{d}^{(0,1)}(e_1,e_2)_- \right)
\end{eqnarray*}
with
$$
(e_1,e_2)_-  = \frac{1}{2}(i_{X_2} \alpha_1 -i_{X_1} \alpha_2
 + f_2 g_1 - f_1 g_2 ).
$$
Dirac structures for the Courant-Jacobi algebroid $({\mathcal
E}^1(M), (0,1)+(0,0))$ will be called {\em Dirac-Jacobi structures}.

\vspace{3mm}

Let us now ``twist" the bracket $\lcf \cdot,\cdot \rcf$ on
$\Gamma({\mathcal E}^1(M))$ with a section $(\varphi, \omega)$ of
$\bigwedge^3(T^*M\times \R)$ by setting
$$
[e_1,e_2]_{(\varphi, \omega)}= \lcf e_1,e_2 \rcf + (\varphi,
\omega)((X_1,f_1) , (X_2,f_2), \cdot).
$$
\begin{proposition}\label{p4.2}
The pair $({\mathcal E}^1(M), (0,1)+(0,0))$ equipped with the
bracket $[\cdot,\cdot]_{(\varphi, \omega)}$ on $\Gamma({\mathcal
E}^1(M))$, the canonical bilinear symmetric form $(\cdot,\cdot)_+$
on the bundle\footnote{$(e_1,e_2)_+  = \frac{1}{2}(i_{X_2} \alpha_1
+i_{X_1} \alpha_2
 + f_2 g_1 + f_1 g_2)$} and the bundle map $\rho=\pi+0$,
is a Courant-Jacobi  algebroid over $M$ if and only if
$\mathrm{d}^{(0,1)}(\varphi, \omega)=0$.

We denote this new Courant-Jacobi algebroid by $({\mathcal
E}^1(M)_{(d \omega,\omega)}, (0,1)+(0,0))$ or simply by $({\mathcal
E}^1(M)_{\omega},(0,1)+(0,0))$.
\end{proposition}
\begin{proof}
We know that $({\mathcal E}^1(M), (0,1)+(0,0))$ equipped with $(\lcf
\cdot,\cdot \rcf, \rho, (\cdot,\cdot)_+)$ is a Courant-Jacobi
algebroid \cite{jj}. Hence, we only have to check the effect of
adding the term $(\varphi, \omega)((X_1,f_1) , (X_2,f_2), \cdot)$ to
the bracket $\lcf \cdot,\cdot \rcf$ on $\Gamma ({\mathcal E}^1(M))$.
Let us set $\theta =(0,1)+(0,0)$. Then, for any
$e_1=(X_1,f_1)+(\alpha_1,g_1), e_2=(X_2,f_2)+(\alpha_2,g_2) \in
\Gamma ({\mathcal E}^1(M))$,  we compute
\begin{eqnarray*}
\rho^{\theta} ([e_1,e_2]_{(\varphi, \omega)}) & = & \rho^{\theta}
(\lcf e_1,e_2 \rcf )+ \rho^{\theta}
((\varphi,\omega)((X_1,f_1),(X_2,f_2), \cdot))  \\
& = & [\rho^{\theta} (e_1), \rho^{\theta} (e_2)],
\end{eqnarray*}
and $ii)$ of Definition \ref{def2} holds. Moreover, for any
$e=(X,f)+(\alpha,g) \in \Gamma ({\mathcal E}^1(M))$, condition
$iii)$ holds if and only if
\begin{eqnarray*}
\lefteqn{((\varphi,\omega)((X,f),(X_1,f_1),
\cdot),(X_2,f_2)+(\alpha_2,g_2))_+ } \\
& +& ((X_1,f_1)+(\alpha_1,g_1),(\varphi,\omega)((X,f),(X_2,f_2),
\cdot))_+   = 0,
\end{eqnarray*}
that is, if and only if
\begin{eqnarray*}
\lefteqn{((\varphi(X,X_1, \cdot) + \omega(f X_1-f_1 X, \cdot ),
\omega(X, X_1)),(X_2,f_2)+(\alpha_2,g_2))_+}  \\ & +&
((X_1,f_1)+(\alpha_1,g_1), (\varphi(X,X_2, \cdot)+ \omega(f X_2-f_2
X, \cdot ), \omega(X, X_2))_+ = 0,
\end{eqnarray*}
which can be proved by a simple computation. Finally, by a long but
straightforward computation, we obtain
\begin{eqnarray*}
[[e_1,e_2]_{(\varphi, \omega)},e_3]_{(\varphi, \omega)} + c.p. & =
& \mathrm{d}^{(0,1)}(T_{(\varphi,\omega)}(e_1,e_2,e_3)) \\
& - & (\mathrm{d}^{(0,1)}(\varphi,
\omega))((X_1,f_1),(X_2,f_2),(X_3,f_3),\cdot)
\end{eqnarray*}
with $T_{(\varphi,\omega)}(e_1,e_2,e_3)=
\frac{1}{3}([e_1,e_2]_{(\varphi, \omega)},e_3)_+ + c.p.$. Thus,
condition $i)$ of Definition \ref{def2} holds if and only if
$\mathrm{d}^{(0,1)}(\varphi, \omega)=(0,0)$ and the proof is
complete.
\end{proof}

\begin{definition}
A Dirac sub-bundle $L$ for the Courant-Jacobi algebroid $({\mathcal
E}^1(M)_\omega, (0,1)+ (0,0))$ over $M$ is called an $\omega${\em
-Dirac-Jacobi structure} or a {\em twisted Dirac-Jacobi structure}.
\end{definition}

Obviously, if $L$ is a twisted Dirac-Jacobi structure, then
$(L,[\cdot,\cdot]_{(d\omega, \omega)}|_L, \rho|_L)$ is a Lie
algebroid over $M$.

\vspace{.2cm}

The next result enables us to characterize twisted Jacobi manifolds
in terms of twisted Dirac-Jacobi structures. Hereafter, in order to
simplify the notation, we will denote the bracket $[\cdot,
\cdot]_{(d\omega, \omega)}$ by $[\cdot, \cdot]_{ \omega}$, whenever
is clear to which bracket we refer to.

\begin{proposition}
Let $ \omega$ be a $2$-form on $M$ and $(\Lambda,E)$ a section of
$\bigwedge^2(TM \times \R)$. Then, {\rm graph}$(\Lambda,E)^\#$ is a
$\omega$-Dirac-Jacobi structure if and only if
$$
[(\Lambda, E), (\Lambda, E)]^{(0,1)}= 2 (\Lambda,E)^{\#}
(d\omega,\omega).
$$
\end{proposition}
\begin{proof}
For any $(\alpha,f), (\beta,g) \in \Gamma(T^*M\times \R)$, we have
\begin{eqnarray*}
\lefteqn{[(\Lambda,E)^\# (\alpha,f) + (\alpha,f), (\Lambda,E)^\#
(\beta,g) + (\beta,g)]_{ \omega}  }  \\
&  &  =[(\Lambda,E)^\# (\alpha,f),(\Lambda,E)^\# (\beta,g)] +
\{(\alpha,f),(\beta,g)\}   \\ & & + (d\omega,\omega)((\Lambda,E)^\#
(\alpha,f),  (\Lambda,E)^\# (\beta,g), \cdot),
\end{eqnarray*}
where $\{ \cdot, \cdot \}$ is the bracket (\ref{1.5}). So, {\rm
graph}$(\Lambda,E)^\#$ is closed under the bracket $[\cdot,\cdot
]_{\omega}$ if and only if
\begin{eqnarray} \label{4.4}
\lefteqn{[(\Lambda,E)^\# (\alpha,f),(\Lambda,E)^\# (\beta,g)] =}
\nonumber \\ & & = (\Lambda,E)^\# (\{(\alpha,f),(\beta,g)\} +
(d\omega,\omega)((\Lambda,E)^\# (\alpha,f),  (\Lambda,E)^\#
(\beta,g), \cdot)).
\end{eqnarray}
But (\ref{4.4}) is equivalent to $[(\Lambda, E), (\Lambda,
E)]^{(0,1)}= 2 (\Lambda,E)^{\#}(d\omega,\omega)$ (see, e.g.
\cite{ks2}).
\end{proof}

\begin{corol} \label{c4.5}
The triple $(M,(\Lambda,E), \omega)$ is a twisted Jacobi manifold if
and only if {\rm graph}$(\Lambda,E)^\#$ is a $\omega$-Dirac-Jacobi
structure.
\end{corol}

Let $(\eta,\gamma)$ be a section of $\bigwedge^2(T^*M \times \R)$.
We denote by $(\eta,\gamma)^\flat: TM \times \R \to T^*M \times \R$
the associated vector bundle morphism that induces on the spaces of
sections a map, that we also denote by $(\eta,\gamma)^\flat$, which
is given, for any $(X,f) \in \Gamma( TM \times \R)$, by
$$
(\eta,\gamma)^\flat(X,f)=(i_X \eta + f \gamma, -i_X \gamma).
$$
\begin{proposition}
Let $(\eta,\gamma)$ be a section of $\bigwedge^2(T^*M \times \R)$.
Then, {\rm graph}$(\eta,\gamma)^{\flat}$ is a $\omega$-Dirac-Jacobi
structure if and only if $\mathrm{d}^{(0,1)}(\eta, \gamma)+ (d
\omega, \omega)=(0,0)$.
\end{proposition}
\begin{proof}
We start by remarking that
$$
\mathrm{d}^{(0,1)}(\eta, \gamma)+ (d \omega, \omega)=(0,0) \,
\Leftrightarrow \, \eta=d \gamma - \omega.
$$
The vector bundle ${\rm graph}(\eta,\gamma)^{\flat}$ over $M$, whose
space of sections is given by
\begin{equation*}
\Gamma({\rm graph}(\eta,\gamma)^{\flat})= \{(X,f)+ (i_X \eta + f
\gamma, -i_X \gamma)\,| \,(X,f) \in \Gamma(TM\times \R)\},
\end{equation*}
is a maximally isotropic sub-bundle of ${\mathcal E}^1(M)$ with
respect to the symmetric bilinear form $(\cdot, \cdot)_+$. Now, let
$e_i=(X_i,f_i)+(i_{X_{i}} \eta + f_i \gamma, -i_{X_i}\gamma)$,
$i=1,2$, be two sections of {\rm graph}$(\eta,\gamma)^{\flat}$.
Then,
\begin{eqnarray*}
[e_1, e_2]_{\omega}& = & \left( [X_1, X_2], X_1(f_2)-X_2(f_1) \right) \\
& & + {\mathcal L}^{(0,1)}_{(X_1,f_1)}(i_{X_{2}} \eta + f_2 \gamma,
-i_{X_2}\gamma) - i_{(X_2,f_2)} \mathrm{d}^{(0,1)}(i_{X_{1}} \eta +
f_1
\gamma, -i_{X_1}\gamma)  \\
& &  +(d \omega, \omega)((X_1,f_1),(X_2,f_2), \cdot)
\end{eqnarray*}
and $[e_1, e_2]_{\omega} \in \Gamma({\rm
graph}(\eta,\gamma)^{\flat})$ if and only if
\begin{eqnarray} \label{4.6}
\lefteqn{{\mathcal L}^{(0,1)}_{(X_1,f_1)}(i_{X_{2}} \eta + f_2
\gamma, -i_{X_2}\gamma) - i_{(X_2,f_2)} \mathrm{d}^{(0,1)}(i_{X_{1}}
\eta + f_1
\gamma, -i_{X_1}\gamma)} \nonumber \\
& & + (d \omega, \omega)((X_1,f_1),(X_2,f_2), \cdot)  \nonumber \\
&=& (i_{[X_1,X_2]} \eta + (X_1(f_2)- X_2(f_1)) \gamma,
-i_{(X_1,X_2]} \gamma).
\end{eqnarray}
A simple computation shows that (\ref{4.6}) is equivalent to $\eta=d
\gamma - \omega$.
\end{proof}

\vspace{3mm}

Let us now look at some other examples of twisted Dirac-Jacobi
structures. Recall that a sub-bundle $L$ of the vector bundle $TM
\oplus T^*M$ over $M$ is a {\em Dirac structure in the sense of
Courant} \cite{c} if $L$ is maximally isotropic with respect to the
symmetric canonical bilinear form on $TM \oplus T^*M$ and
$\Gamma(L)$ closes under the Courant bracket, which is given,  for
any sections $X+\alpha$, $Y+\beta$ of $TM \oplus T^*M$, by
\begin{equation} \label{cbracket}
[X+ \alpha, Y + \beta ]_C= [X,Y]+ {\mathcal L}_X \beta - {\mathcal
L}_Y \alpha + \frac{1}{2} d(i_Y \alpha-i_X \beta).
\end{equation}

\begin{example}
{\em Let $L$ be a sub-bundle of $TM \oplus T^*M$, $\omega$ a
$2$-form on $M$ and consider the sub-bundle $L_\omega$ of ${\mathcal
E}_{\omega}^1(M)$ whose fiber at a point $x\in M$ is given by
$$
L_{ \omega}(x)= \{ (X, 0)_x +(\alpha - i_X \omega , f)_x \, | \, (X
+ \alpha)_x \in L_x \}.
$$
Then, $L_\omega$ is a $\omega$-Dirac-Jacobi structure if and only if
$L$ is a Dirac structure in the sense of Courant. It is immediate to
verify that $L_\omega$ is maximally isotropic with respect to
symmetric canonical bilinear form on ${\mathcal E}_{\omega}^1(M)$ if
and only if $L$ is maximally isotropic with respect to symmetric
canonical bilinear form on $TM \oplus T^*M$. Moreover, if
$(X,0)+(\alpha-i_X \omega,f)$ and $(Y,0)+ (\beta -i_Y \omega, g)$
are any two sections of $L_\omega$, then
\begin{eqnarray*}
\lefteqn{[(X,0)+(\alpha-i_X \omega,f),(Y,0)+ (\beta -i_Y \omega,
g)]_\omega =([X,Y],0)+} \\
& & + ({\mathcal L}_X \beta - {\mathcal L}_Y \alpha + \frac{1}{2}
d(i_Y \alpha-i_X \beta)- i_{[X,Y]}\omega,X.g-Y.f +\frac{1}{2}d(i_Y
\alpha-i_X \beta)).
\end{eqnarray*}
So, the sections of $L_\omega$ close under the bracket $[\cdot,
\cdot]_\omega$ if and only if the sections of $L$ close under the
Courant bracket on $TM \oplus T^*M$.}
\end{example}

For the next example we need the following definition.

\begin{definition}
A {\em twisted locally conformal presymplectic} structure on a
manifold $M$ is a pair $((\Theta, \vartheta), \omega)$, where
$\Theta$ and $\omega$ are two $2$-forms on $M$ and $\vartheta$ is a
closed $1$-form on $M$ such that
$$
d(\Theta + \omega) + \vartheta \wedge (\Theta + \omega) =0.
$$
If $M$ is even dimensional and $\Theta$ is non-degenerate,
$(M,(\Theta, \vartheta), \omega)$ is a twisted locally conformal
symplectic manifold (cf. Example \ref{exemples-tJ}.2).
\end{definition}

\begin{example}
{\em Let $\Theta$ and $\omega$ be two $2$-forms on a manifold $M$
and $\vartheta$ be a $1$-form on $M$. Consider the sub-bundle
$L_{((\Theta, \vartheta), \omega)}$ of ${\mathcal E}_{\omega}^1(M)$
whose fiber at a point $x\in M$ is given by
\begin{equation}  \label{tlcs}
L_{((\Theta, \vartheta), \omega)}(x)= \{ (X, i_X \vartheta)_x
+(i_X\Theta - f \vartheta, f)_x \, | \, (X,f)_x\in (TM\times \R)_x
\}.
\end{equation}
Then, $L_{((\Theta, \vartheta), \omega)}$ is a twisted Dirac-Jacobi
structure if and only if $((\Theta, \vartheta), \omega)$ is a
twisted locally conformal presymplectic structure on $M$.
Effectively, it is easy to check that $L_{((\Theta, \vartheta),
\omega)}$ is a maximally isotropic sub-bundle of ${\mathcal
E}_{\omega}^1(M)$, with respect to the bilinear symmetric form
$(\cdot,\cdot)_+$. Let $(X, i_X \vartheta) +(i_X \Theta - f
\vartheta, f)$ and $(Y, i_Y \vartheta)+(i_Y \Theta -g \vartheta, g)$
be two sections of $L_{((\Theta, \vartheta), \omega)}$. We compute
\begin{eqnarray*}
\lefteqn{[(X, i_X \vartheta) +(i_X \Theta - f \vartheta, f), (Y, i_Y
\vartheta)+(i_Y
\Theta -g \vartheta, g)]_{\omega}} \\
&=& ([X,Y], i_{[X,Y]}\vartheta + d\vartheta(X,Y))+\\
& & + (i_{[X,Y]}\Theta - g i_X d \vartheta + f i_Y d \vartheta +
d\Theta(X,Y,\cdot)+ d \omega (X,Y, \cdot)\\
& & +(\vartheta \wedge \Theta)(X,Y, \cdot) + (\vartheta \wedge
\omega)(X,Y, \cdot)\\
& & - \{X.g -Y.f-(i_X \vartheta)g+ (i_Y
\vartheta)f+ \Theta(X,Y)+ \omega(X,Y)\}\vartheta, \\
& &X.g -Y.f-(i_X \vartheta)g+ (i_Y \vartheta)f+ \Theta(X,Y)+
\omega(X,Y));
\end{eqnarray*}
so, the space $\Gamma(L_{((\Theta, \vartheta), \omega)})$ is closed
under the bracket $[\cdot, \cdot]_\omega$ if and only if
$d\vartheta=0$ and $d(\Theta + \omega) + \vartheta \wedge (\Theta +
\omega) =0$.}
\end{example}

For the next example, we recall that if $(M, \Lambda, \varphi)$ is a
twisted Poisson manifold, then $(T^*M, \{ \cdot, \cdot \}^{\varphi},
\Lambda^\#)$ is a Lie algebroid over $M$, where the Lie bracket $\{
\cdot, \cdot \}^{\varphi}$ is defined, for all $1$-forms $\alpha$
and $\beta$ on $M$, by
$$\{ \alpha, \beta \}^{\varphi}={\mathcal L}_{\Lambda^\#(\alpha)}
\beta - {\mathcal L}_{\Lambda^\#(\beta)} \alpha - d(\Lambda(\alpha,
\beta)) + \varphi (\Lambda^\#(\alpha),\Lambda^\#(\beta), \cdot).$$

\begin{example}
{\em Let $\Lambda$ be a bivector filed on $M$, $Z$ a vector field on
$M$ and $\omega$ a $2$-form on $M$. We denote by $L_{(\Lambda, Z,
\omega)}$ the sub-bundle of ${\mathcal E}_{\omega}^1(M)$ whose fiber
at a point $x\in M$ is given by}
$$
L_{(\Lambda, Z, \omega)} (x)= \{(\Lambda^\#(\alpha) - f Z, f)_x
+(\alpha, i_Z \alpha)_x \, | \, (\alpha,f)_x\in (T^*M\times \R)_x
\}.
$$
{\em Then, $L_{(\Lambda, Z, \omega)}$ is a twisted Dirac-Jacobi
structure if and only if $(\Lambda, Z, \omega)$ defines an
homogeneous twisted exact Poisson structure on $M$ (cf. Definition
\ref{def.tehp}). An easy computation shows that $L_{(\Lambda, Z,
\omega)}$ is a maximally isotropic sub-bundle of ${\mathcal
E}^1(M)$, with respect to symmetric bilinear form $(\cdot,
\cdot)_+$. Let $(\Lambda^\#(\alpha) - f Z, f) +(\alpha, i_Z \alpha)$
and $(\Lambda^\#(\beta) - g Z, g) +(\beta, i_Z \beta)$ be two
sections of $L_{(\Lambda, Z, \omega)}$. Then, if $(\Lambda, Z,
\omega)$ is a twisted exact homogeneous Poisson structure, we
compute}
\begin{eqnarray*}
\lefteqn{[(\Lambda^\#(\alpha) - f Z, f) +(\alpha, i_Z
\alpha),(\Lambda^\#(\beta) - g Z, g) +(\beta, i_Z
\beta)]_{\omega}}\\
& = & \left( \Lambda^\# \left( {\mathcal L}_{\Lambda^\#(\alpha)}
\beta - {\mathcal L}_{\Lambda^\#(\beta)} \alpha- d ( \Lambda(\alpha,
\beta))+ g ({\mathcal L}_Z \alpha - \alpha)- f({\mathcal L}_Z \beta
- \beta)\right. \right.\\
& & + \left. d \omega(\Lambda^\#(\alpha),\Lambda^\#(\beta), \cdot)
\right) - \left((\Lambda^\#(\alpha)).g- (\Lambda^\#(\beta)).f +
g(Z.f)-f (Z.g)\right)Z\, ,\\
& &  \left. (\Lambda^\#(\alpha)).g- (\Lambda^\#(\beta)).f +
g(Z.f)-f (Z.g) \right) + \\
& + & \left({\mathcal L}_{\Lambda^\#(\alpha)} \beta - {\mathcal
L}_{\Lambda^\#(\beta)} \alpha- d ( \Lambda(\alpha, \beta))+ g
({\mathcal L}_Z \alpha - \alpha)- f({\mathcal L}_Z \beta - \beta)
\right. \\
& & + d
\omega(\Lambda^\#(\alpha),\Lambda^\#(\beta), \cdot)\\
& & \underbrace{-d\omega(\Lambda^\#(\alpha), g Z, \cdot)+
d\omega(\Lambda^\#(\beta), f Z, \cdot) - \omega(g\Lambda^\#(\alpha),
\cdot) + \omega(f\Lambda^\#(\beta), \cdot)}_{=0} \, , \\
& & i_Z \left({\mathcal L}_{\Lambda^\#(\alpha)} \beta - {\mathcal
L}_{\Lambda^\#(\beta)} \alpha- d ( \Lambda(\alpha, \beta))+ g
({\mathcal L}_Z \alpha - \alpha)- f({\mathcal L}_Z \beta -
\beta) \right. \\
& & \left. \left. +d \omega(\Lambda^\#(\alpha),\Lambda^\#(\beta),
\cdot) \right) \right)
\end{eqnarray*}
{\em and we conclude that the space of sections of $L_{(\Lambda, Z,
\omega)}$ is closed under the bracket $[\cdot, \cdot ]_\omega$.
Thus, $L_{(\Lambda, Z, \omega)}$ is a $\omega$-Dirac-Jacobi
structure. A similar computation shows that, conversely, if
$L_{(\Lambda, Z, \omega)}$ is a $\omega$-Dirac-Jacobi structure,
then the triple $(\Lambda, Z, \omega)$ defines an homogeneous
twisted exact Poisson structure on $M$.}
\end{example}

Let $\varphi$ be a closed $3$-form on $M$ and $L$ a sub-bundle of
$TM \oplus T^*M$. We recall that $L$ is called a $\varphi${\em
-Dirac structure} (in the sense of Courant) \cite{sw} if it is
maximally isotropic with respect to the canonical bilinear symmetric
form on $TM \oplus T^*M$, and its space of sections is closed under
the bracket $[\cdot, \cdot]_{C_\varphi}$ which is given, for any
sections $X+\alpha$ and $Y+\beta$ of $TM \oplus T^*M$, by
$$
[X+\alpha,Y+\beta]_{C_\varphi}=[X+\alpha,Y+\beta]_C + \varphi(X,Y,
\cdot),
$$
where $[ \cdot, \cdot]_C$ is the Courant bracket given by
(\ref{cbracket}). In \cite{im2},  the existence of a correspondence
between Dirac-Jacobi structures $L \subset {\mathcal E}^1(M)$ and
Dirac structures $\tilde{L}\subset T(M\times \R) \oplus T^*(M\times
\R)$ in the sense of Courant was proved (see also \cite{fj}). For
twisted Dirac-Jacobi structures we can establish the following.

\begin{proposition}
Let $L$ be a sub-bundle of ${\mathcal E}_{\omega}^1(M)$. Then,
\begin{enumerate}
\item the sub-bundle $\tilde{L}_{\omega}\subset T(M\times \R)
\oplus T^*(M\times \R)$ given by
$$
\tilde{L}_{\omega}= \{(X+f \frac{\partial}{\partial t})+ e^t(\alpha
+i_X \omega +g dt) \, | \, (X,f)+(\alpha,g) \in L \}
$$
is a Dirac structure (in the sense of Courant) if and only if $L$ is
an $\omega$-Dirac-Jacobi structure; \item the sub-bundle
$\tilde{L}\subset T(M\times \R) \oplus T^*(M\times \R)$ given by
$$
\tilde{L}= \{(X+f \frac{\partial}{\partial t})+ e^t(\alpha +g dt) \,
| \, (X,f)+(\alpha,g) \in L \}
$$
is a $d(e^t \omega)$-Dirac structure (in the sense of Courant)if and
only if $L$ is an $\omega$-Dirac-Jacobi structure.
\end{enumerate}
\end{proposition}
\begin{proof}
A simple computation proves that each one of the sub-bundles
$\tilde{L}_{\omega}$ and $\tilde{L}$ of $T(M\times \R) \oplus
T^*(M\times \R)$ is maximally isotropic with respect to the
canonical symmetric bilinear form in $T(M\times \R) \oplus
T^*(M\times \R)$ if and only if $L \subset {\mathcal
E}_{\omega}^1(M)$ is maximally isotropic with respect to the
canonical symmetric bilinear form in ${\mathcal E}_{\omega}^1(M)$.
To complete the proof of the first assertion, we take two sections
$(X_i + f_i \displaystyle{\frac{\partial}{\partial t}})+e^t(\alpha_i
+ i_{X_i}\omega + g_i dt)$, $i=1,2$, of $\tilde{L}_{\omega}$. Then,
denoting by $[\cdot , \cdot]_C$ the Courant bracket on $\Gamma(T (M
\times \R) \oplus T^* (M \times \R))$, we compute
\begin{eqnarray*}
\lefteqn{[(X_1 + f_1 \displaystyle{\frac{\partial}{\partial
t}})+e^t(\alpha_1 + i_{X_1}\omega + g_1 dt), (X_2 + f_2
\displaystyle{\frac{\partial}{\partial t}})+e^t(\alpha_2 +
i_{X_2}\omega + g_2 dt)]_C} \\
 & = & \left( [X_1, X_2]+ (X_1 .f_2
-X_2.f_1)\displaystyle{\frac{\partial}{\partial t}} \right) + e^t
\left( {\mathcal L}_{X_1} \alpha_2 - {\mathcal L}_{X_2} \alpha_1
\right. \\
& & + \frac{1}{2}d(i_{X_2}\alpha_1 - i_{X_1}\alpha_2) +f_1 \alpha_2
- f_2 \alpha_1 + \frac{1}{2}(g_2 df_1-g_1df_2-f_1dg_2+f_2dg_1) \\
& & \left. + d \omega(X_1,X_2, \cdot) + \omega(f_1X_2-f_2X_1, \cdot)
+
i_{[X_1,X_2]} \omega \right)\\
& & + e^t \left( X_1.g_2- X_2.g_1 + \frac{1}{2}(i_{X_2}\alpha_1-
i_{X_1}\alpha_2 -f_2g_1+f_1g_2) +\omega(X_1,X_2)\right) dt
\end{eqnarray*}
and, since
\begin{eqnarray*}
\lefteqn{([X_1, X_2], X_1 .f_2 -X_2.f_1)}\\
& & + \left({\mathcal L}_{X_1} \alpha_2 - {\mathcal L}_{X_2}
\alpha_1 + \frac{1}{2}d(i_{X_2}\alpha_1 - i_{X_1}\alpha_2)+f_1
\alpha_2 - f_2 \alpha_1 \right.\\
& & +  \frac{1}{2}(g_2 df_1-g_1df_2-f_1dg_2+f_2dg_1) + d
\omega(X_1,X_2, \cdot) +
\omega(f_1X_2-f_2X_1, \cdot),\\
& & \left. X_1.g_2- X_2.g_1 + \frac{1}{2}(i_{X_2}\alpha_1-
i_{X_1}\alpha_2 -f_2g_1+f_1g_2) +\omega(X_1,X_2) \right)\\
& = &[(X_1,f_1)+(\alpha_1,g_1), (X_2,f_2)+(\alpha_2,g_2)]_\omega,
\end{eqnarray*}
we conclude that the bracket $[\cdot, \cdot]_C$ closes in
$\Gamma(\tilde{L}_{\omega})$ if and only if the bracket $[\cdot,
\cdot]_\omega$ closes in $\Gamma(L)$. The proof of the second
assertion is very similar and we omit it.
\end{proof}

\section{Gauge transformations}
As in \cite{sw}, in the case of (twisted) Dirac structures for
Courant algebroids, we may define gauge transformations for
Dirac-Jacobi sub-bundles. Given a section $(\eta, \gamma)$ of
$\bigwedge^2(T^*M \times \R)$, let us consider the vector bundle map
$$\tau_{(\eta, \gamma)} : {\mathcal E}^1(M) \to {\mathcal E}^1(M)$$
\noindent that induces on the  spaces of sections a map, that we
also denote by  $\tau_{(\eta, \gamma)}$, which is defined, for any
$(X,f)+(\alpha,g) \in \Gamma({\mathcal E}^1(M))$, by
$$
\tau_{(\eta, \gamma)}((X,f)+(\alpha,g))=(X,f)+(\alpha,g)+(\eta,
\gamma)^{\flat}(X,f).
$$
$\tau_{(\eta, \gamma)}$ is called a {\em gauge transformation}
associated with $(\eta, \gamma)$. Let us also consider the
Courant-Jacobi algebroids  $({\mathcal E}^1(M)_{(d \omega,\omega)},
(0,1)+(0,0))$ and  $({\mathcal E}^1(M)_{(d
\omega,\omega)-\mathrm{d}^{(0,1)}(\eta, \gamma)}, (0,1)+(0,0))$.
Given a $(d \omega,\omega)$-Dirac-Jacobi structure $L$, its image by
$\tau_{(\eta, \gamma)}$ is the vector sub-bundle of ${\mathcal
E}^1(M)$,
$$
\tau_{(\eta, \gamma)}(L)= \{(X,f)+(\alpha,g)+(\eta,
\gamma)^{\flat}(X,f)\,| \, (X,f)+(\alpha,g)\in L \}.
$$

\begin{proposition}
Let $L$ be a $(d \omega,\omega)$-Dirac-Jacobi structure. Then, for
any $(\eta, \gamma) \in \Gamma(\bigwedge^2(T^*M \times \R))$,
$\tau_{(\eta, \gamma)}(L)$ is a $((d
\omega,\omega)-\mathrm{d}^{(0,1)}(\eta,\gamma))$-Dirac-Jacobi
structure. Moreover,
$$
\tau_{(\eta, \gamma)}|_L: (L,[\cdot,\cdot]_{(d \omega,
\omega)}|_L,\rho|_L) \to (\tau_{(\eta, \gamma)}(L),[\cdot,\cdot]_{(d
\omega, \omega)-\mathrm{d}^{(0,1)}(\eta, \gamma)}|_{\tau_{(\eta,
\gamma)}(L)},\rho|_{\tau_{(\eta, \gamma)}(L)} )
$$
is an isomorphism of Lie algebroids over the identity, with $\rho=
\pi + 0$.
\end{proposition}
\begin{proof}
Let $e_1=(X_1,f_1)+(\alpha_1,g_1)$ and
$e_2=(X_2,f_2)+(\alpha_2,g_2)$ be any two sections of $L$. Then,
\begin{eqnarray*}
(\tau_{(\eta, \gamma)}(e_1), \tau_{(\eta, \gamma)}(e_2))_+ & = &
\underbrace{(e_1,e_2)_+}_{=0} + \frac{1}{2}(\langle
(X_1,f_1),(\eta, \gamma)^{\flat}(X_2,f_2) \rangle \nonumber \\
& &  + \langle(X_2,f_2),(\eta, \gamma)^{\flat}(X_1,f_1)\rangle ) =0,
\end{eqnarray*}
and $\tau_{(\eta, \gamma)}(L)$ is a maximally isotropic sub-bundle
of ${\mathcal E}^1(M)$. On the other hand,
\begin{eqnarray} \label{4.19}
\lefteqn{[\tau_{(\eta, \gamma)}(e_1), \tau_{(\eta,
\gamma)}(e_2)]_{(d \omega, \omega)-\mathrm{d}^{(0,1)}(\eta,
\gamma)}= \lcf
e_1, e_2 \rcf} \nonumber \\
& &+(d(\omega- \eta), \omega- \eta+d \gamma)((X_1,f_1), (X_2,f_2),
\cdot)\nonumber \\
& &+ {\mathcal L}^{(0,1)}_{(X_1,f_1)}(i_{X_2} \eta +f_2 \gamma,
-i_{X_2} \gamma)-{\mathcal L}^{(0,1)}_{(X_2,f_2)}(i_{X_1} \eta +f_1
\gamma, -i_{X_1} \gamma)\nonumber \\ & & +
\mathrm{d}^{(0,1)}(i_{X_2}(i_{X_1}\eta) + f_1(i_{X_2}\gamma)- f_2(i_{X_1}\gamma))\nonumber \\
&=& \lcf e_1, e_2 \rcf +((d\omega, \omega)-(d \eta, \eta-d
\gamma))((X_1,f_1), (X_2,f_2),
\cdot)\nonumber \\
& & + (i_{[X_1,X_2]} \eta +i_{X_2}(i_{X_1}d \eta)+ (X_1.f_2)\gamma
+f_2({\mathcal L}_{X_1} \gamma) +f_1 (i_{X_2} \eta)\nonumber \\& &
-f_2 (i_{X_1} \eta)-(X_2.f_1)\gamma -f_1(i_{X_2} d\gamma)- f_2 d
(i_{X_1}\gamma),\nonumber \\
& & \eta(X_1,X_2) - i_{X_1}(i_{X_2}\gamma)+ i_{X_2}(i_{X_1}\gamma)) \nonumber \\
&=& \lcf e_1, e_2 \rcf +(d\omega, \omega)((X_1,f_1), (X_2,f_2),
\cdot)+ i_{[(X_1,f_1), (X_2,f_2)]^{(0,1)}}(\eta, \gamma)\nonumber \\
& = &\tau_{(\eta, \gamma)}([e_1,e_2]_{(d\omega,\omega)}),
\end{eqnarray}
which means that $\Gamma(\tau_{(\eta, \gamma)}(L))$ closes under the
bracket $[\cdot , \cdot ]_{(d \omega,
\omega)-\mathrm{d}^{(0,1)}(\eta, \gamma)}$ and we conclude that
$\tau_{(\eta, \gamma)}(L)$ is a $((d
\omega,\omega)-\mathrm{d}^{(0,1)}(\eta,\gamma))$-Dirac-Jacobi
structure.

Moreover, with $\rho= \pi + 0$, we have
\begin{equation}  \label{4.20}
\rho(\tau_{(\eta, \gamma)}((X,f)+(\alpha,g)))=
\rho((X,f)+(\alpha,g)),
\end{equation}
for any section $(X,f)+(\alpha,g)$ of $L$. From (\ref{4.19}) and
(\ref{4.20}), we deduce that $\tau_{(\eta, \gamma)}|_L$ is an
isomorphism of Lie algebroids over the identity.
\end{proof}

\vspace{3mm}

The twisted Dirac-Jacobi structures $L$ and $\tau_{(\eta,
\gamma)}(L)$ are said to be {\em gauge-equivalent}.

\begin{corol} \label{corol4.8}
Let $L$ be a $(d \omega,\omega)$-Dirac-Jacobi structure and $(\eta,
\gamma) \in \Gamma(\bigwedge^2(T^*M \times \R))$.
\begin{itemize}
\item[\bf i)] If $\mathrm{d}^{(0,1)}(\eta, \gamma)=(0,0)$, then
$\tau_{(\eta, \gamma)}(L)$ is also a $(d
\omega,\omega)$-Dirac-Jacobi structure. \item[\bf ii)] If
$\mathrm{d}^{(0,1)}(\eta, \gamma)=(d \omega, \omega)$, then
$\tau_{(\eta, \gamma)}(L)$ is a Dirac-Jacobi structure.
\end{itemize}
\end{corol}

Let us denote by {\em Dir}$_\omega$ the set of all
$\omega$-Dirac-Jacobi structures and consider the additive group
\begin{equation*}
{\mathcal F}= \{ (\eta ,\gamma) \in \Gamma(\bigwedge^2(T^*M \times
\R))\, | \, \mathrm{d}^{(0,1)}(\eta ,\gamma)=0 \}.
\end{equation*}
Corollary \ref{corol4.8} i) means that ${\mathcal F}$ acts on {\em
Dir}$_\omega$ with the action,
$$
{\mathcal F} \times  Dir_{\omega} \rightarrow  Dir{_\omega}, \qquad
((\eta ,\gamma),L) \mapsto \tau_{(\eta, \gamma)}(L),
$$
and two elements of {\em Dir}$_\omega$ are gauge equivalent if they
lie in the same orbit of the action.

\section{The Jacobi algebroid associated to a twisted Jacobi manifold}
In this section we will show that we can associate a Jacobi
algebroid to each twisted Jacobi manifold.

\begin{proposition} \label{p5.1}
Let $(M,(\Lambda,E),\omega)$ be a twisted Jacobi manifold. Then,
$(T^*M\times \R, \{\cdot,\cdot \}^{\omega},\pi \circ
(\Lambda,E)^\#)$ is a Lie algebroid over $M$, where $\{\cdot,\cdot
\}^{\omega}$ is the bracket on $\Gamma(T^*M\times \R)$ given, for
all $(\alpha,f),(\beta,g) \in \Gamma(T^*M\times \R)$, by
\begin{equation}\label{5.1}
\{ (\alpha,f),(\beta,g) \} ^{\omega} = \{(\alpha,f),(\beta,g)\} + (d
\omega,\omega)((\Lambda,E)^\# (\alpha,f), (\Lambda,E)^\# (\beta,g),
\cdot),
\end{equation}
$\{\cdot,\cdot\}$ being the bracket (\ref{1.5}).
\end{proposition}
\begin{proof}
Let $(M,(\Lambda,E),\omega)$ be a twisted Jacobi manifold. From
Corollary \ref{c4.5}, we know that {\rm graph}$(\Lambda,E)^\#$ is a
twisted Dirac-Jacobi sub-bundle of $\mathcal{E}_{\omega}^1(M)$,
hence it is a Lie algebroid over $M$ with the following bracket on
the space of its sections,
\begin{eqnarray} \label{5.2}
\lefteqn{[(\Lambda,E)^\# (\alpha,f) + (\alpha,f), (\Lambda,E)^\#
(\beta,g) + (\beta,g)]_{(d \omega, \omega)}  } \nonumber \\
&  &  =[(\Lambda,E)^\# (\alpha,f),(\Lambda,E)^\# (\beta,g)]
\nonumber \\ & & + \{(\alpha,f),(\beta,g)\}  + (d
\omega,\omega)((\Lambda,E)^\# (\alpha,f),  (\Lambda,E)^\# (\beta,g),
\cdot).
\end{eqnarray}
Since the bracket (\ref{5.2}) splits in the sum $\Gamma(TM \times
\R) \oplus \Gamma(T^*M \times \R)$, then its projection $\{ \cdot,
\cdot \} ^{\omega}$ over $\Gamma(T^*M \times \R)$ is a Lie bracket.
Moreover, for any $h \in \C$,
$$
\{(\alpha,f),h(\beta,g)\}^{\omega} = h \{
(\alpha,f),(\beta,g)\}^{\omega} + (((\pi \circ (\Lambda,E)^\#)
(\alpha,f))\cdot h )(\beta,g).
$$
So, $(\{\cdot, \cdot \}^{\omega},\pi \circ (\Lambda,E)^\#)$ endows
$T^*M\times \R$ with a Lie algebroid structure.
\end{proof}

\begin{corol}
Let $(M,(\Lambda,E), \omega)$ be a twisted Jacobi manifold. Then,
for any $f,g \in \C$,
$$
\{ \mathrm{d}^{(0,1)}f,\mathrm{d}^{(0,1)}g \} ^{\omega} =
\mathrm{d}^{(0,1)} \{f,g \} +(d \omega, \omega)((\Lambda,E)^\#
(\mathrm{d}^{(0,1)}f), (\Lambda,E)^\# (\mathrm{d}^{(0,1)}g), \cdot).
$$
\end{corol}
\begin{proof}
It is an immediate consequence of Proposition \ref{p5.1}, taking
into account that, for any $f,g \in \C$, $\mathrm{d}^{(0,1)} \{f,g
\}= \{ \mathrm{d}^{(0,1)}f,\mathrm{d}^{(0,1)}g \}$, with the bracket
on the left hand-side given by (\ref{br-M}) and the bracket on the
right hand-side given by (\ref{1.5}).
\end{proof}

\vspace{.3cm}

The differential operator $d_*^{\omega}$ defined on
$\Gamma(\bigwedge (TM\times \R))$ by the Lie algebroid structure
$(\{\cdot,\cdot\}^{\omega}, \pi \circ (\Lambda,E)^\#)$ on $T^*M
\times \R$ is given,
\begin{itemize}
\item for any $f \in \C$, by
\begin{equation} \label{5.10}
d_*^{\omega}f= d_*f= -(\Lambda,E)^\# (df,0);
\end{equation}
\item for any $(X,f) \in \Gamma(TM \times \R)$, by
\begin{equation} \label{5.11}
d_*^{\omega}(X,f)=d_*(X,f) +( (\Lambda,E)^{\#} \otimes 1) (d \omega,
\omega)(X,f),
\end{equation}
\end{itemize}
where $d_*$ denotes the operator given by (\ref{1.6}) and
$(\Lambda,E)^{\#} \otimes 1$ is defined adapting (\ref{2.3}) in the
obvious way.

\begin{proposition} \label{p5.2}
Let $(M,(\Lambda,E), \omega)$ be a twisted Jacobi manifold. The
section $(-E,0)$ of $TM \times \R$ is a $1$-cocycle for the Lie
algebroid $(T^*M\times \R, \{\cdot,\cdot \}^{\omega},\pi \circ
(\Lambda,E)^\#)$ over $M$.
\end{proposition}
\begin{proof}
It suffices to prove that $d_*^{\omega}(-E,0)=(0,0)$. Let
$(\alpha,f),(\beta,g)$ be any sections of $T^*M \times \R$. Then,
\begin{eqnarray*}
\lefteqn{ d_*^{\omega}(-E,0)((\alpha,f),(\beta,g))}
\nonumber \\
& = & d_*(-E,0)((\alpha,f),(\beta,g)) +( (\Lambda,E)^{\#} \otimes 1)
(d \omega, \omega)(-E,0)((\alpha,f),(\beta,g))
\nonumber \\
& = & [E, \Lambda](\alpha, \beta)-(d \omega,
\omega)((\Lambda,E)^\#(\alpha,f), (\Lambda,E)^\#(\beta,g), (-E,0))
\nonumber \\
& = & (( \Lambda^{\#} \otimes 1)(d \omega)(E)- (( \Lambda^{\#}
\otimes 1)(\omega)(E)) \wedge
E)(\alpha, \beta) \nonumber \\
& & +d \omega(\Lambda^\#(\alpha), \Lambda^\#(\beta), E)- (i_E
\alpha) \omega(\Lambda^\#(\beta),E)+(i_E \beta)
\omega(\Lambda^\#(\alpha),E)
\nonumber \\
& = & 0,
\end{eqnarray*}
and so, $d_*^{\omega}(-E,0)=(0,0)$.
\end{proof}

\vspace{.2cm}

From Propositions \ref{p5.1} and \ref{p5.2}, we deduce that the
twisted Jacobi structure $((\Lambda,E), \omega)$ on $M$ defines a
Jacobi algebroid structure on $T^*M \times \R$. Moreover we have,
from (\ref{1.10}), (\ref{5.10}) and (\ref{5.11}), that
\begin{itemize}
\item for any $f \in \C$,
\begin{equation} \label{7.5}
(d_*^{ \omega})^{(-E,0)}f= -(\Lambda,E)^\#(df, f);
\end{equation}
\item for any $(X,f) \in \Gamma(TM \times \R)$,
\begin{equation}  \label{7.6}
(d_*^{ \omega})^{(-E,0)}(X,f)=[(\Lambda,E),
(X,f)]^{(0,1)}+((\Lambda,E)^{\#} \otimes 1) (d \omega, \omega)(X,f).
\end{equation}
\end{itemize}

The Lie algebra homomorphism, from $\C$ to $\Gamma(TM)$, expressed
by equation (\ref{1.4}) in the case where $M$ is a Jacobi manifold,
fails in the case of twisted Jacobi manifolds, as shown in the next
proposition.

\begin{proposition}
Let $(M,(\Lambda,E),  \omega)$ be a twisted Jacobi manifold. Then,
for any $f,g \in \C$,
\begin{equation}  \label{5.3}
[X_f, X_g]= X_{ \{f,g \} }+ (\pi \circ (\Lambda,E)^\#)( (d \omega,
\omega)((\Lambda,E)^\#(df,f), (\Lambda,E)^\#(dg,g), \cdot)).
\end{equation}
\end{proposition}
\begin{proof}
From (\ref{4.4}) we have, with $(\Lambda,E)^\# (df,f)=(X_f, -E\cdot
f)$,
\begin{eqnarray*}
\lefteqn{[(X_f,-E\cdot f), (X_g,-E\cdot g)]=(X_{ \{f,g \}}, -E\cdot
\{f,g \})}  \\ & & + (\Lambda,E)^\#( (d \omega,
\omega)((\Lambda,E)^\#(df,f), (\Lambda,E)^\#(dg,g), \cdot)).
\end{eqnarray*}
The projection over the first factor gives (\ref{5.3}).
\end{proof}

\section{Quasi-Jacobi bialgebroids and their doubles}

The notion of {\em quasi-Lie bialgebroid} was introduced in
\cite{ro}. It is a structure on a pair $(A,A^*)$ of vector bundles,
in duality, over a differentiable manifold $M$ that is defined by a
Lie algebroid structure on $A^*$, a skew-symmetric bracket on the
space of smooth sections of $A$ and a bundle map $a : A \to TM$,
satisfying some compatibility conditions. These conditions are
expressed in terms of a section of $\bigwedge^3A^*$, which turns to
be an obstruction to the Lie bialgebroid structure on $(A,A^*)$. As
in the case of a Lie bialgebroid, the double $A\oplus A^*$ of a
quasi-Lie bialgebroid $(A,A^*)$ is endowed with  a Courant algebroid
structure \cite{ro, ks3}.

In this section, in order to adapt the previous notion to the Jacobi
framework, we introduce the concept of \emph{quasi-Jacobi
bialgebroid} and we prove that its double is endowed with a
Courant-Jacobi algebroid structure \cite{jj, gm2}.

\begin{definition} \label{d6.1}
A {\em quasi-Jacobi bialgebroid} structure on a pair $(A,A^*)$ of
dual vector bundles over a differentiable manifold $M$ consists of:
\begin{itemize}
\item a Lie algebroid structure $([\cdot,\cdot]_*, a_*)$ on $A^*$
with a $1$-cocycle $W$; \item a bundle map $a:A \to TM$; \item a
skew-symmetric operation $[\cdot, \cdot]$ on $\Gamma(A)$; \item a
section $\phi \in \Gamma(A^*)$; \item a section $\varphi \in
\Gamma(\bigwedge^3A^*)$;
\end{itemize}
satisfying, for all $X,Y,Z \in \Gamma(A)$ and $f \in
 \C$, the following properties:
 \begin{enumerate}
 \item [1)]
 $[X,fY]=f[X,Y]+(a(X)f)Y$;
 \item [2)]
 $a([X,Y])=[a(X),a(Y)]-a_* \varphi(X,Y,\cdot)$;
\item [3)] $[[X,Y],Z]+ c.p. = - d_*^W(\varphi(X,Y,Z)) -
((i_{\varphi(X,Y, \cdot)}d_*^W Z)+ c.p.)$, where $d_*^W$ is the
modified differential operator on $\Gamma(\bigwedge A)$ defined by
the Lie algebroid structure of $A^*$ and the 1-cocycle $W$; \item
[4)] $d \phi - \varphi(W,\cdot,\cdot)=0$, where $d$ is the
quasi-differential operator on $\Gamma(\bigwedge A^*)$ determined by
the structure $([\cdot, \cdot],a)$ on $A$; \item [5)] $d^\phi
\varphi=0$, \noindent where $d^\phi$ is given, for any $\beta \in
\Gamma(\bigwedge^k A^*)$, by $d^\phi(\beta)= d \beta + \phi \wedge
\beta$; \item [6)] $d_{*}^{W} [P,Q]^{\phi}= [d_{*}^{W} P,Q]^{\phi} +
(-1)^{p+1}[P, d_{*}^{W}Q]^{\phi}$, with $P \in \Gamma
(\bigwedge^{p}A)$ and $Q \in \Gamma (\bigwedge A)$.
\end{enumerate}
We will denote the quasi-Jacobi bialgebroid by $((A,
\phi),(A^*,W),\varphi)$.
\end{definition}

Let $((A, \phi),(A^*,W),\varphi)$ be a quasi-Jacobi bialgebroid over
$M$, $\mathcal{L}^\phi$ and $\mathcal{L}_*^W$ the quasi-Lie
derivative and the Lie derivative operators defined, respectively,
by $d^\phi$ and $d_*^W$ as in (\ref{1.11}), $a^\phi$ and $a_*^W$ the
deformed anchor maps according to (\ref{1.12}). On the Whitney sum
bundle $A\oplus A^*$ we consider the two nondegenerate canonical
bilinear forms $(\cdot,\cdot)_\pm$ and, on the space $\Gamma(A\oplus
A^*)\cong \Gamma(A)\oplus \Gamma(A^*)$ we define the bracket
$[\cdot,\cdot]_\varphi$ by setting, for any $e_1=X_1+\alpha_1,
e_2=X_2+\alpha_2 \in \Gamma(A\oplus A^*)$,
\begin{equation}
[e_1,e_2]_{\varphi}=[X_1+ \alpha_1, X_2+ \alpha_2]_{\varphi}= \lcf
X_1+ \alpha_1, X_2+ \alpha_2 \rcf +\varphi (X_1,X_2, \cdot),
\end{equation}
where $\lcf \cdot, \cdot \rcf$ is the bracket (\ref{4.3}).

\begin{theorem}\label{th7}
Let $((A, \phi),(A^*,W),\varphi)$ be a quasi-Jacobi bialgebroid over
$M$. The vector bundle $A\oplus A^*$ over $M$ endowed with
$([\cdot,\cdot]_{\varphi}, (\cdot,\cdot)_+,
\rho^{\theta},\mathcal{D}^\theta)$, where $\theta = \phi +W \in
\Gamma(A^*\oplus A)$, $\rho^\theta = a^\phi + a_*^W$ and
$\mathcal{D}^\theta = (d_*^W + d^\phi)\vert_{\C}$, is a
Courant-Jacobi algebroid over $M$.
\end{theorem}

For establishing the above theorem, we need the results of the
following lemmas. Let $((A,\phi),(A^*,W),\varphi)$ be a quasi-Jacobi
bialgebroid over $M$.

\begin{lemma}\label{lem7.1}
For any $P \in \Gamma(\bigwedge^k A)$, $X,Y \in \Gamma(A)$, $\alpha
\in \Gamma(A^*)$ and $f \in \C$,
\begin{itemize}
\item [i)] $d_*^W [X,Y]= [d_*^W X, Y]+ [X, d_*^W Y]$; \item [ii)]
${\mathcal L}_{* \phi}^{W}P + \, {\mathcal L}_W^{\phi } P=0$;
\item [iii)] $\langle \phi , W \rangle=0 \quad and \quad
a(W)+a_*(\phi)=0$; \item[iv)] ${\mathcal L}_{* \phi} X + [W,X]=0$;
\item[v)] $[d_*^W f, X] + {\mathcal L}^W_{* d^\phi f}X=0 \quad and
\quad [d^{\phi}f,\alpha]_*^W + \mathcal{L}_{d_*^Wf}^\phi \alpha =
0$.
\end{itemize}
\end{lemma}
\begin{proof}
The proof is based on the facts that $d_*^W$ (resp. $d^\phi$) is a
derivation of $[\cdot,\cdot]^\phi$ (resp. $[\cdot,\cdot]_*^W$) and
it is similar to the case of a Jacobi bialgebroid (see
\cite{im1,jj}).
\end{proof}

\vspace{2mm}

On the space $\C$ we define the internal composition law
$\{\cdot,\cdot\}$ by setting, for any $f,g \in \C$,
\begin{equation}  \label{88.1}
\{f,g \}=\langle d^\phi f, d_*^W g\rangle.
\end{equation}

\begin{lemma}  \label{l8.1}
For any $f,g \in \C$,
\begin{equation}
[d_*^W f, d_*^W g]= d_*^W (\{g,f \}).
\end{equation}
\end{lemma}
\begin{proof}
From the skew-symmetry of the bracket $[\cdot,\cdot]$ on
$\Gamma(A)$, from Lemma \ref{lem7.1} $v)$ and because $(d_*^W)^2=0$,
$$
[d_*^W f, d_*^W g]  =  - [d_*^W g, d_*^W f]= {\mathcal L}^W_{*
d^\phi g} (d_*^W f)  =  d_*^W (\langle d^\phi g, d_*^W f\rangle)=
d_*^W ( \{g,f \}).
$$
\end{proof}

\begin{lemma} \label{l8.2}
The bracket (\ref{88.1}) is a first-order differential operator on
the second argument and it is skew-symmetric.
\end{lemma}
\begin{proof}
In fact, for any $f,g,h \in \C$,
\begin{equation}  \label{8.20}
\{f, gh \}=g \{f, h \}+ h \{f, g \} - gh \{f, 1 \}
\end{equation}
because
$$
d_*^W(gh)= g \,d_*^W h + h\, d_*^W g -gh W.
$$

In order to establish the skew-symmetry of (\ref{8.1}), we will
prove that, for any $f \in \C$,
\begin{equation}\label{8.0}
\{f,f \}=0.
\end{equation}
Since $(A^*, [\cdot,\cdot]_*, a_*, W)$ is a Lie algebroid over $M$
with a $1$-cocycle, the homomorphism of $\C$-modules $a_*^W :
\Gamma(A^*) \to \Gamma(TM \times \R)$ given by (\ref{1.12}), induces
a Lie algebroid homomorphism over the identity between the Lie
algebroids with 1-cocycles $(A^*, [\cdot,\cdot]_*, a_*, W)$ and
$(TM\times \R, [\cdot, \cdot], \pi, (0,1))$. Hence, for any $f\in
\C$ \footnote{In this section, in order to avoid confusion with the
quasi-differential $d$ of $A$, we will denote by $\delta f$ the
usual de Rham differential of $f\in \C$.},
\begin{equation}\label{8.2}
(a_*^W)^*(0,1)=W, \quad (a_*^W)^* (\delta f,f)=d_*^W f \quad
\mathrm{and} \quad (a_*^W)^* (\delta f,0)= d_* f,
\end{equation}
where $(a_*^W)^* : \Gamma (T^*M\times \R) \to \Gamma(A)$ denotes the
transpose of $a_*^W$. On the other hand, since the
quasi-differential operator $d$ on $\Gamma(A^*)$ is defined by $a :
\Gamma(A) \to \Gamma(TM)$ and by the bracket $[\cdot,\cdot]$ on
$\Gamma(A)$, we can easily prove that
\begin{equation}\label{8.3}
(a^\phi)^* (\delta f,0) = a^*(\delta f) = df \quad \mathrm{and}
\quad (a^\phi)^* (\delta f,f)= d^\phi f,
\end{equation}
where $(a^\phi)^* :\Gamma (T^*M\times \R) \to \Gamma(A^*)$ denotes
the transpose of $a^\phi$. So,
\begin{eqnarray}\label{8.4}
\{f,g\} & = & \langle d^\phi f, d_*^Wg \rangle
\stackrel{(\ref{8.2}),(\ref{8.3})}{=} \langle (a^\phi)^* (\delta
f,f), (a_*^W)^* (\delta g,g)\rangle \nonumber \\
& = & \langle (\delta f,f), a^\phi \circ (a_*^W)^*(\delta
g,g)\rangle.
\end{eqnarray}
When $g=1$, (\ref{8.4}) gives
\begin{equation}\label{8.5}
\{f, 1 \}= \langle (\delta f,f), a^\phi \circ (a_*^W )^* (0,1)
\rangle \stackrel{(\ref{8.2})}{=} \langle (\delta f,f), a^\phi
(W)\rangle = -\langle \delta f, a_*(\phi)\rangle,
\end{equation}
where the last equality follows from Lemma \ref{lem7.1} $iii)$. On
the other hand,
\begin{eqnarray} \label{8.6}
\{1, f \} &= & \langle (0,1),a^\phi \circ (a_*^W )^* (\delta
f,f)\rangle = \langle (0,1),a^\phi(
d_*^W f)\rangle \nonumber \\
& =& \langle \phi, d_*f\rangle=\langle \phi, a_*^*(\delta f)\rangle.
\end{eqnarray}
From (\ref{8.5}) and (\ref{8.6}), we get
\begin{equation} \label{8.7}
\{f, 1 \}=-\{1, f \}.
\end{equation}
Using Lemma \ref{lem7.1} $iii)$, (\ref{8.2}) and (\ref{8.3}), we can
write
\begin{equation} \label{8.8}
\{f,f \}=\langle (\delta f,0),a^\phi \circ (a_*^W )^* (\delta
f,0)\rangle.
\end{equation}
From Lemma \ref{l8.1} we have,
\begin{equation} \label{8.9}
d_*^W( \{f , f \})= [d_*^W f, d_*^W f]=0.
\end{equation}
In particular, for $f^2$,
\begin{equation} \label{8.10}
d_*^W( \{f^2 , f^2 \})=0
\end{equation}
and
\begin{eqnarray*}
0 = d_*^W( \{f^2 , f^2 \}) &\stackrel{(\ref{8.8})}{=}& d_*^W
(\langle (\delta f^2,0), a^\phi \circ (a_*^W )^* (\delta f^2,0)\rangle  \\
& = & 4f^2 d_*^W( \{f , f \})+4 \{f , f \}d_*^W f^2 -4f^2\{f,f\}W \\
&\stackrel{(\ref{8.9}),(\ref{1.10})}{=}& 4 \{f , f \} d_* f^2.
\end{eqnarray*}
So, for any $f \in \C$,
\begin{equation} \label{8.11}
 \{f , f \}d_*f^2=0.
\end{equation}
Then,
\begin{eqnarray}\label{8.12}
0& \stackrel{(\ref{8.11}),(\ref{1.10})}{=}& \langle d^\phi 1,\{f , f
\}d_*^W f^2-f^2\{f,f\}W\rangle  \nonumber \\ & = &\{f , f \}\{1 ,
f^2 \}-f^2\{f,f\}\langle \phi,W\rangle
\nonumber \\
&\stackrel{(\ref{8.20})}{=}& 2f\{f , f \}\{1 , f \} - f^2 \{f , f
\} \underbrace{\{1 , 1 \}}_{=0} \nonumber \\
& = & 2f\{f , f \}\{1 , f \}
\end{eqnarray}
and
\begin{eqnarray*}
0& \stackrel{(\ref{8.11}),(\ref{1.10})}{=} & \langle d^\phi f,\{f ,
f \}d_*^W f^2 - f^2\{f,f\}W\rangle \\& = &\{f , f \}
\{f , f^2 \} -f^2\{f,f\}\{f,1\}\\
&\stackrel{(\ref{8.20}), (\ref{8.7})}{=}& 2 f\{f , f \}^2 +2f^2\{f ,
f \}\{1 , f \} \\
&\stackrel{(\ref{8.12})}{=}& 2f\{f , f \}^2,
\end{eqnarray*}
whence we deduce that (\ref{8.0}) holds.
\end{proof}

\begin{remark}\label{rem7}
{\rm From the skew-symmetry of (\ref{88.1}) and the fact that it is
first order differential operator on the second argument, we
conclude that it is first order differential operator on each
argument.}
\end{remark}

\begin{lemma}\label{lem7.2}
For any $f\in \C$, $X\in \Gamma(A)$ and $\alpha \in \Gamma(A^*)$,
\begin{itemize}
\item[i)] $(a \circ d_*^W + a_*\circ d^\phi)f = 0$; \item[ii)]
$[a(X), a_*(\alpha)]=a_*({\mathcal L}_X^{\phi} \alpha)- a({\mathcal
L}^W_{* \alpha} X)+a(d_*^W \langle \alpha, X \rangle)$.
\end{itemize}
\end{lemma}
\begin{proof}
For $i)$ we have that, for any $g\in \C$,
\begin{eqnarray*}
\langle (a^\phi \circ d_*^W + a_*^W\circ d^\phi)f, (\delta
g,g)\rangle & = & \langle d_*^Wf, (a^\phi)^*(\delta g,g)\rangle +
\langle d^\phi f, (a_*^W)^*(\delta g,g)\rangle \\
& \stackrel{(\ref{8.3}),(\ref{8.2})}{=}&\langle d_*^Wf, d^\phi
g\rangle + \langle d^\phi f, d_*^Wg\rangle \\
& \stackrel{(\ref{88.1})}{=}& \{g,f\} + \{f,g\} = 0,
\end{eqnarray*}
because $\{\cdot,\cdot\}$ is skew-symmetric. So, $(a^\phi \circ
d_*^W + a_*^W\circ d^\phi)f =0$. But,
$$
(a^\phi \circ d_*^W + a_*^W\circ d^\phi)f = (a \circ d_*^W +
a_*\circ d^\phi)f + \langle \phi, d_*^Wf\rangle + \langle W, d^\phi
f\rangle
$$
and
$$
\langle \phi, d_*^Wf\rangle + \langle W, d^\phi f\rangle
\stackrel{(\ref{8.2}),(\ref{8.3})}{=} \langle a_*^W(\phi)+a^\phi(W),
(\delta f,f)\rangle = 0,
$$
where the last equality follows from Lemma \ref{lem7.1} $iii)$.
Consequently, for any $f\in \C$,
$$
(a \circ d_*^W + a_*\circ d^\phi)f = 0.
$$
The proof of $ii)$ is similar to the case of a Jacobi bialgebroid
(see \cite{im1,jj}).
\end{proof}

\begin{lemma}\label{lem-Lie-der}
Let $((A,\phi),(A^*,W),\varphi)$ be a quasi-Jacobi bialgebroid over
$M$. Then, the quasi-Lie derivative operator $\mathcal{L}^\phi$
associated to the quasi-differential operator $d^\phi$ on
$\Gamma(\bigwedge A^*)$ satisfies the following property: For any
$X,Y,V_1, \ldots, V_p \in \Gamma(A)$ and any $\eta \in
\Gamma(\bigwedge^p A^*)$,
\begin{eqnarray}\label{Lie-der}
{\mathcal L}^{\phi}_{[X,Y]} \eta (V_1, \ldots, V_p) & = & ({\mathcal
L}^{\phi}_X \circ {\mathcal L}^{\phi}_Y - {\mathcal L}^{\phi}_Y
\circ {\mathcal L}^{\phi}_X)\eta (V_1, \ldots, V_p)
\nonumber \\
&  & +\sum_{i=1}^{p} (-1)^i \eta([[X,Y],V_i]+ c.p.\, ,V_1, \ldots,
\hat{V_i}, \ldots, V_p) \nonumber \\
&  & - a_*^W(\varphi(X,Y,\cdot))(\eta(V_1,\ldots,V_p)).
\end{eqnarray}
\end{lemma}
\begin{proof}
We prove the above formula by a simple, but long, computation,
taking into account the condition $4)$ of Definition \ref{d6.1} of a
quasi-Jacobi bialgebroid.
\end{proof}

\vspace{3mm}

Now, we will prove Theorem \ref{th7}.

\vspace{3mm}

\begin{proof-of} {\bf{Theorem \ref{th7}.}} We have to check that the conditions $i)-iv)$ of Definition \ref{def2}
hold. In order to establish condition $ii)$, we use the results of
Lemma \ref{lem7.2} and the conditions $2)$ and $4)$ of Definition
\ref{d6.1} of a quasi-Jacobi bialgebroid. We obtain that, for any
two sections $e_1=X_1 + \alpha_1$, $e_2=X_2 + \alpha_2$ of $A \oplus
A^*$,
$$
\rho^\theta ([e_1, e_2]_{\varphi}) =[\rho^\theta (e_1), \rho^\theta
(e_2)].
$$
For condition $iii)$ we have that, for all $e,e_1,e_2 \in \Gamma(A
\oplus A^*)$, $e=X+\alpha$, $e_1=X_1 + \alpha_1$, $e_2=X_2 +
\alpha_2$,
\begin{eqnarray*}
& & ([e,e_1]_{\varphi} + {\mathcal D}^\theta(e,e_1)_+,e_2)_+ +
(e_1,[e,e_2]_\varphi+{\mathcal D}^\theta(e,e_2)_+)_+
\\
& = & (\lcf e,e_1 \rcf +{\mathcal D}^\theta(e,e_1)_+,e_2)_+
+\frac{1}{2}\varphi(X,X_1,X_2) \\
& + & (e_1, \lcf e,e_2 \rcf +{\mathcal D}^{\theta}(e,e_2)_+)_+
+\frac{1}{2}\varphi(X,X_2,X_1) \\
& = & (\lcf e,e_1 \rcf +{\mathcal D}^\theta(e,e_1)_+,e_2)_+ + (e_1,
\lcf e,e_2 \rcf +{\mathcal D}^{\theta}(e,e_2)_+)_+.
\end{eqnarray*}
But, by doing the same computations as in Proposition 4.1 of
\cite{jj}, we establish the equality
$$
(\lcf e,e_1 \rcf +{\mathcal D}^\theta(e,e_1)_+,e_2)_+ + (e_1, \lcf
e,e_2 \rcf +{\mathcal D}^{\theta}(e,e_2)_+)_+ =
\rho^{\theta}(e)(e_1,e_2)_+.
$$
Hence, we conclude
$$
\rho^{\theta}(e)(e_1,e_2)_+ = ([e,e_1]_{\varphi} + {\mathcal
D}^\theta(e,e_1)_+,e_2)_+ + (e_1,[e,e_2]_\varphi+{\mathcal
D}^\theta(e,e_2)_+)_+.
$$
The condition $iv)$ can be easily proved as follows. For any $f,g
\in \C$,
\begin{eqnarray*}
({\mathcal D}^{\theta}f, {\mathcal D}^{\theta}g)_+ &=&(d_*^W
f+d^{\phi}f, d_*^W g+d^{\phi}g)_+  = \frac{1}{2}( \langle
d^{\phi}g,d_*^W f \rangle +  \langle d^{\phi}f,d_*^W g
\rangle )\\
& = & \frac{1}{2}(\{g,f \}+ \{f,g \} )=0,
\end{eqnarray*}
where $\{\cdot, \cdot \}$ is the  bracket (\ref{88.1}) which, by
Lemma \ref{l8.2}, is skew-symmetric. Finally, it remains to
establish condition $i)$ of Definition \ref{def2}, i.e., for any
$e_1,e_2,e_3 \in \Gamma(A\oplus A^*)$, $e_i=X_i+\alpha_i$,
$i=1,2,3$,
\begin{equation}\label{Jac-id}
[[e_1,e_2]_\varphi,e_3]_\varphi + [[e_2,e_3]_\varphi,e_1]_\varphi +
[[e_3,e_1]_\varphi,e_2]_\varphi = \mathcal{D}^\theta
T_\varphi(e_1,e_2,e_3),
\end{equation}
where
$T_{\varphi}(e_1,e_2,e_3)=\displaystyle{\frac{1}{3}(([e_1,e_2]_{\varphi}
,e_3)_+ + c.p.)}$. Since the proof involves a very long computation,
we only give a short schedule.

First, we note that, if
$T(e_1,e_2,e_3)=\displaystyle{\frac{1}{3}((\lcf e_1,e_2\rcf,e_3)_+ +
c.p.)}$, then
\begin{equation}\label{T}
T_{\varphi}(e_1,e_2,e_3)= T(e_1,e_2,e_3) +
\frac{1}{2}\varphi(X_1,X_2,X_3).
\end{equation}
Let us set
\begin{equation}\label{Y-B}
[[e_1,e_2]_\varphi,e_3]_\varphi + [[e_2,e_3]_\varphi,e_1]_\varphi +
[[e_3,e_1]_\varphi,e_2]_\varphi = Y+\beta,
\end{equation}
where $Y$ and $\beta$ denote the components of
$[[e_1,e_2]_\varphi,e_3]_\varphi +c.p.$ on $\Gamma(A)$ and
$\Gamma(A^*)$, respectively. We have
\begin{eqnarray*}
[[e_1,e_2]_\varphi,e_3]_\varphi + c.p. & = & [\lcf e_1,e_2\rcf
+\varphi(X_1,X_2,\cdot),e_3]_\varphi + c.p. \\
& = & ([\lcf e_1,e_2\rcf ,e_3]_\varphi +
[\varphi(X_1,X_2,\cdot),e_3]_\varphi)+ c.p. \\
&=& (\lcf \lcf e_1,e_2\rcf ,e_3\rcf +\varphi(\widetilde{\lcf
e_1,e_2\rcf},X_3,\cdot)+\lcf \varphi(X_1,X_2,\cdot),e_3\rcf)+c.p.,
\end{eqnarray*}
where $\widetilde{\lcf e_i,e_j\rcf}$, $i,j=1,2,3$, denotes the part
of $\lcf e_i,e_j\rcf$ that belongs to $\Gamma(A)$. Hence,
$$
Y = (\widetilde{\lcf \lcf e_1,e_2\rcf ,e_3\rcf}+ \widetilde{\lcf
\varphi(X_1,X_2,\cdot),e_3\rcf})+c.p.
$$
Taking into account condition $3)$ of Definition \ref{d6.1}, the
fact that $(A^*,[\cdot,\cdot]_*,a_*)$ is a Lie algebroid over $M$,
so
$$
\mathcal{L}^W_{* [\alpha_i,\alpha_j]_*^W} = \mathcal{L}^W_{*
\alpha_i}\circ \mathcal{L}^W_{* \alpha_j} - \mathcal{L}^W_{*
\alpha_j}\circ \mathcal{L}^W_{* \alpha_i}, \quad \mathrm{for} \quad
i,j=1,2,3,
$$
and also (\ref{T}), we obtain, after a long computation,
\begin{equation}\label{Y}
Y=d_*^W(T_\varphi(e_1,e_2,e_3)).
\end{equation}
Similarly, for $\beta$ we have
$$
\beta = (\widehat{\lcf \lcf e_1,e_2\rcf ,e_3\rcf} +
\varphi(\widetilde{\lcf e_1,e_2\rcf},X_3,\cdot) + \widehat{\lcf
\varphi(X_1,X_2,\cdot),e_3\rcf})+c.p.,
$$
where $\widehat{\lcf \lcf e_i,e_j\rcf ,e_k\rcf}$ (resp.
$\widehat{\lcf \varphi(X_i,X_j,\cdot),e_k\rcf}$), $i,j,k=1,2,3$,
denotes the component of $\lcf \lcf e_i,e_j\rcf ,e_k\rcf$ (resp.
$\lcf \varphi(X_i,X_j,\cdot),e_k\rcf$) that is section of $A^*$. We
repeat the computations developed in Proposition 4.1 of \cite{jj}
for the calculation of the corresponding $\beta$ and we take into
account the conditions $3)$ and $5)$ of Definition \ref{d6.1}, the
fact that $(A^*,[\cdot,\cdot]_*,a_*)$ is a Lie algebroid over $M$,
so $[[\alpha_1,\alpha_2]_*^W]_*^W + c.p. =0$, the result of Lemma
\ref{lem-Lie-der} and (\ref{T}). After a long calculation we get
\begin{equation}\label{B}
\beta = d^\phi(T_\varphi(e_1,e_2,e_3)).
\end{equation}
From (\ref{Y-B}), (\ref{Y}) and (\ref{B}) we conclude that
(\ref{Jac-id}) holds.
\end{proof-of}
\begin{remark}
{\rm When $\varphi=0$, the quasi-Jacobi bialgebroid is a Jacobi
bialgebroid and we obtain Proposition 4.1 of \cite {jj}.}
\end{remark}

\section{The quasi-Jacobi bialgebroid of a twisted Jacobi
manifold}

Let $(M,(\Lambda,E),\omega)$ be a twisted Jacobi manifold. We
consider the following skew-symmetric bracket on the space of
sections of the vector bundle $TM \times \R$ over $M$, given, for
all $(X,f), (Y,g) \in \Gamma(TM \times \R)$, by
\begin{equation} \label{7.1}
[(X,f), (Y,g)]'= [(X,f), (Y,g)]-(\Lambda,E)^\#((d \omega,
\omega)((X,f), (Y,g), \cdot)),
\end{equation}
where $[\cdot, \cdot]$ is the bracket (\ref{1.7}), and we define an
operator $\textrm{d}'$, acting on the space of sections of the
exterior algebra $\bigwedge (T^*M \times \R)$ as a graduate
differential operator, by setting,
\begin{itemize}
\item on $f \in \C$,
$$
\textrm{d}'f= \textrm{d}f=(df,0);
$$
\item on sections $(\alpha,f)$ of $T^*M \times \R$,
$$
\textrm{d}'(\alpha,f)=\textrm{d}(\alpha,f) -(d\omega,
\omega)((\Lambda,E)^\#(\alpha,f),\cdot,\cdot).
$$
\end{itemize}
Then, we extend $\textrm{d}'$, by linearity, to the algebra
$(\Gamma(\bigwedge (T^*M \times \R)),\wedge)$. The operator
$\textrm{d}'$ coincides with the one determined by the structure
$([\cdot, \cdot]', \pi)$ on $TM \times \R$.

Now, we use the section $(0,1) \in \Gamma(T^*M \times \R)$ to modify
the bracket $[\cdot , \cdot ]'$ on $\Gamma(TM \times \R)$, according
to formula (\ref{defor}), and also the operator $\textrm{d}'$. The
new bracket will be denoted by $[\cdot,\cdot]'^{(0,1)}$ and the
resulting operator $\textrm{d}'^{(0,1)}$ is defined as follows:
\begin{itemize}
\item on $f \in \C$,
$$
\textrm{d}'^{(0,1)} f=\mathrm{d}^{(0,1)}f=(df,f);
$$
\item on sections $(\alpha,f)$ of $T^*M \times \R$,
$$
\textrm{d}'^{(0,1)}(\alpha,f)=\mathrm{d}^{(0,1)}(\alpha,f)
-(d\omega, \omega)((\Lambda,E)^\#(\alpha,f),\cdot,\cdot).
$$
\end{itemize}

Let us extend the bracket $[\cdot,\cdot]'^{(0,1)}$ on $\Gamma(TM
\times \R)$ to the whole algebra  $(\Gamma(\bigwedge(TM \times \R)),
\wedge)$, as in the case of a Jacobi algebroid.
 In particular, if $(X,f) \in \Gamma(TM \times
\R)$ and $(C,Y) \in \Gamma(\bigwedge^2(TM \times \R))$, we have
\begin{eqnarray} \label{8.1}
[(C,Y),(X,f)]'^{(0,1)} & = &[(C,Y),(X,f)]^{(0,1)} \nonumber \\ &- &
\left(((\Lambda, E)^\# \otimes(C,Y)^\# + (C,Y)^\# \otimes (\Lambda,
E)^\#)\otimes 1
\right) (d\omega, \omega)(X,f),\nonumber \\
\end{eqnarray}
where the second term of the right hand-side of (\ref{8.1}) is the
section of $\bigwedge^2(TM \times \R)$ given, for any $(\alpha,g),
(\beta,h)\in \Gamma(T^*M \times \R)$, by
\begin{eqnarray*}
& &  \left(((\Lambda, E)^\# \otimes(C,Y)^\# + (C,Y)^\# \otimes
(\Lambda, E)^\#)\otimes 1 \right) (d\omega, \omega)(X,f)((\alpha,
g), (\beta, h))= \\
&  & \quad =(d\omega, \omega)((\Lambda, E)^\#(\alpha,
g),(C,Y)^\#(\beta, h),(X,f)) \\
& & \qquad +(d\omega, \omega)((C, Y)^\#(\alpha,
g),(\Lambda,E)^\#(\beta, h),(X,f)).
\end{eqnarray*}

\begin{lemma}
Let $(M, (\Lambda,E),\omega)$ be a twisted Jacobi manifold. Then,
for any $(X,f) \in \Gamma(TM \times \R)$, we have
\begin{equation} \label{88.2}
(d_*^{ \omega})^{(-E,0)}(X,f)=[(\Lambda,E), (X,f)]'^{(0,1)} -
((\Lambda,E)^\# \otimes 1)(d\omega, \omega)(X,f).
\end{equation}
\end{lemma}
\begin{proof}
It is a direct consequence of (\ref{8.1}), (\ref{defor}) and
(\ref{7.6}).
\end{proof}

\vspace{.3cm}

We remark that if $\omega$ is a $2$-form on $M$ such that
$(\Lambda,E)^\#(d\omega,\omega)=(0,0)$, i.e. when the twisted Jacobi
manifold is just a Jacobi manifold, we recover the well-known
relation \cite{im1}, $d_* ^{(-E,0)}(X,f)=[(\Lambda,E),
(X,f)]^{(0,1)}$.

\vspace{.3cm}

In the next theorem, which is the main result of this section, we
show that one can associate a quasi-Jacobi bialgebroid to each
twisted Jacobi manifold.

\begin{theorem}  \label{p7.1}
Let $(M, (\Lambda,E),\omega)$ be a twisted Jacobi manifold and
$(T^*M \times \bkR, \{ \cdot,\cdot \}^{\omega}, \pi \circ
(\Lambda,E)^\#)$ its associated Lie algebroid. Consider the vector
bundle $TM \times \R$ equipped with the bracket (\ref{7.1}) on the
space of its sections, the operator ${\rm d}'$ and the projection
$\pi: TM \times \R \to TM$. Then, $((TM \times \R, (0,1)),(T^*M
\times \R, (-E,0)), (d\omega,\omega))$ is a quasi-Jacobi bialgebroid
over $M$.
\end{theorem}
\begin{proof}
We have to check that all conditions of Definition \ref{d6.1} are
satisfied. According to Proposition \ref{p5.2}, the section $(-E,0)$
of  $TM \times \R$ is a $1$-cocycle for the Lie algebroid $(T^*M
\times \bkR, \{ \cdot,\cdot \}^{\omega}, \pi \circ (\Lambda,E)^\#)$.

Let  $(X,f)$ and $(Y,g)$ be any two sections of  $TM \times \R$ and
$h \in \C$. Then,
$$
[(X,f), h(Y,g)]'= h [(X,f), (Y,g)]' + (\pi(X,f))(h) (Y,g),
$$
which means that condition 1) of Definition \ref{d6.1} holds. We
also have
$$
\pi([(X,f),(Y,g)]')=[X,Y]- (\pi \circ (\Lambda,E)^\#)((d \omega,
\omega)((X,f), (Y,g), \cdot)),
$$
which is $2)$ of Definition \ref{d6.1}. Moreover,
$$
{\rm d}'(0,1)= -(d \omega, \omega)((\Lambda,E)^\#(0,1), \cdot,
\cdot)=(d \omega, \omega)((-E,0),\cdot, \cdot)
$$
and so 4) is also satisfied. The skew-symmetry of the morphism
$(\Lambda,E)^\#$, allows us to conclude that
$$
{\rm d}'^{(0,1)}(d \omega, \omega)=(0,0),
$$
which is condition 5) of Definition \ref{d6.1}.

Let us now consider the sections $(X_1,f_1)$, $(X_2,f_2)$ and
$(X_3,f_3)$ of $TM  \times \R$. Then,
\begin{eqnarray} \label{7.20}
\lefteqn{[[(X_1,f_1), (X_2,f_2)]', (X_3,f_3)]' + c.p.  =
\left([[(X_1,f_1),
(X_2,f_2)], (X_3,f_3)] \right.}\nonumber \\
& & \left.- (\Lambda,E)^\#((d \omega, \omega)([(X_1,f_1),
(X_2,f_2)],(X_3,f_3), \cdot)) \right. \nonumber \\
& & \left.- [(\Lambda,E)^\#((d \omega, \omega)((X_1,f_1), (X_2,f_2),
\cdot)),(X_3,f_3)] \right. \nonumber \\
& & \left.+ (\Lambda,E)^\#((d \omega, \omega)((\Lambda,E)^\#((d
\omega, \omega)((X_1,f_1), (X_2,f_2), \cdot)),(X_3,f_3), \cdot))
\right) + c.p. \nonumber \\
\end{eqnarray}
First, we remark that, since $[\cdot, \cdot]$ is a Lie bracket on
$\Gamma(TM \times \R)$,
$$
[[(X_1,f_1), (X_2,f_2)], (X_3,f_3)]+ c.p.=(0,0).
$$
Let $(\alpha,g)$ be an arbitrary section of $T^*M \times \R$. Then,
\begin{eqnarray}  \label{7.30}
\lefteqn{\langle(\alpha,g),-(\Lambda,E)^\#((d \omega,
\omega)([(X_1,f_1),
(X_2,f_2)],(X_3,f_3), \cdot))+ c.p.\rangle  }\nonumber \\
& =& (\pi(X_1,f_1)).((d \omega, \omega)((X_2,f_2), (X_3,f_3),
(\Lambda,E)^\#(\alpha,g)))+c.p.\nonumber \\
& & -(\pi(\Lambda,E)^\#(\alpha,g)).((d \omega,
\omega)((X_1,f_1),(X_2,f_2), (X_3,f_3)))\nonumber \\
& & -(d \omega,
\omega)([(X_1,f_1),(\Lambda,E)^\#(\alpha,g)],(X_2,f_2),
(X_3,f_3))-c.p.\nonumber \\
& & + (0,1) \wedge (d \omega, \omega)((X_1,f_1),(X_2,f_2),
(X_3,f_3),(\Lambda,E)^\#(\alpha,g))
\end{eqnarray}
and
\begin{eqnarray} \label{7.31}
\lefteqn{\langle(\alpha,g),- [(\Lambda,E)^\#((d \omega,
\omega)((X_1,f_1),
(X_2,f_2), \cdot)),(X_3,f_3)] +c.p.\rangle }\nonumber \\
&= &- \left( i_{(d \omega,
\omega)((X_1,f_1),(X_2,f_2),\cdot)}(d_{*}^\omega)^{(-E,0)}(X_3,f_3)
\right)(\alpha,g) -c.p.\nonumber \\
& & -(\pi(X_1,f_1)).((d \omega, \omega)((X_2,f_2), (X_3,f_3),
(\Lambda,E)^\#(\alpha,g)))-c.p.\nonumber \\
& & +(d \omega,
\omega)([(X_1,f_1),(\Lambda,E)^\#(\alpha,g)],(X_2,f_2),
(X_3,f_3))+c.p.\nonumber \\
& & -f_1 (d \omega, \omega)((X_2,f_2), (X_3,f_3),(\Lambda,E)^\#(\alpha,g))-c.p.\nonumber \\
& & -(d \omega, \omega)((X_1,f_1), (X_2,f_2),((\Lambda,E)^\# \otimes
1)(d \omega, \omega)(X_3,f_3)(\alpha,g)) -c.p.
\end{eqnarray}
On the other hand,
\begin{eqnarray} \label{7.32}
\lefteqn{\langle(\alpha,g),(\Lambda,E)^\#((d \omega,
\omega)((\Lambda,E)^\#((d \omega, \omega)((X_1,f_1), (X_2,f_2),
\cdot)),(X_3,f_3), \cdot)) + c.p.\rangle}\nonumber \\
& =&(d \omega, \omega)((X_1,f_1), (X_2,f_2),((\Lambda,E)^\# \otimes
1)(d \omega, \omega)(X_3,f_3)(\alpha,g)) +c.p.
\end{eqnarray}
If we add up the terms of (\ref{7.30}), (\ref{7.31}) and
(\ref{7.32}), we obtain
\begin{eqnarray*}
\lefteqn{-(\pi(\Lambda,E)^\#(\alpha,g)).((d \omega,
\omega)((X_1,f_1),(X_2,f_2), (X_3,f_3)))}\nonumber \\
& & + (0, d
\omega)((X_1,f_1),(X_2,f_2), (X_3,f_3),(\Lambda,E)^\#(\alpha,g))\nonumber \\
& & -\left( i_{(d \omega,
\omega)((X_1,f_1),(X_2,f_2),\cdot)}(d_{*}^\omega)^{(-E,0)}(X_3,f_3)
\right)(\alpha,g) -c.p.\nonumber \\
& & -f_1 (d \omega, \omega)((X_2,f_2),
(X_3,f_3),(\Lambda,E)^\#(\alpha,g))-c.p. \nonumber \\
& &= -(d_*^\omega)^{(-E,0)}((d \omega, \omega)((X_1,f_1),(X_2,f_2),
(X_3,f_3)))((\alpha,g))\nonumber \\
& &-\left( i_{(d \omega,
\omega)((X_1,f_1),(X_2,f_2),\cdot)}(d_{*}^\omega)^{(-E,0)}(X_3,f_3)
\right)(\alpha,g) -c.p.
\end{eqnarray*}
and we conclude that
\begin{eqnarray*}
\lefteqn{[[(X_1,f_1), (X_2,f_2)]', (X_3,f_3)]' + c.p. =}
\\
& & -(d_*^\omega)^{(-E,0)}((d \omega, \omega)((X_1,f_1),(X_2,f_2),
(X_3,f_3))) \\
& &-\left( i_{(d \omega,
\omega)((X_1,f_1),(X_2,f_2),\cdot)}(d_{*}^\omega)^{(-E,0)}(X_3,f_3)
 +c.p.\right),
\end{eqnarray*}
which is condition 3) of Definition \ref{d6.1}.

Finally, we must show that, for any $(P,P_0) \in
\Gamma(\bigwedge^p(TM\times  \R))$ and $(Q,Q_0) \in
\Gamma(\bigwedge(TM\times \R))$,
\begin{eqnarray} \label{deriv}
\lefteqn{(d_*^{ \omega})^{(-E,0)}[(P,P_0), (Q,Q_0)]'^{(0,1)}=} \nonumber \\
&=& [(d_*^{ \omega})^{(-E,0)}(P,P_0), (Q,Q_0)]'^{(0,1)} + (-1)^{p+1}
[(P,P_0),(d_*^{ \omega})^{(-E,0)}(Q,Q_0)]'^{(0,1)}.
\end{eqnarray}
As in the case of a Jacobi algebroid \cite{gm1}, it is enough to
prove (\ref{deriv}) in the cases where: i) $(P,P_0)$ and $(Q,Q_0)$
are both functions of $M$; ii) $(P,P_0)$ is a section of $TM \times
\R$ and $(Q,Q_0)$ is a function of $M$; iii) $(P,P_0)$ and $(Q,Q_0)$
are both sections of $TM \times \R$.

\vspace{.2cm}

We remark that, for any $f \in \C$ and $(P,P_0) \in
\Gamma(\bigwedge^p(TM\times  \R))$,
$$
[(P,P_0), f]'^{(0,1)}= [(P,P_0), f]^{(0,1)}.
$$
When $(P,P_0)=(f,0)\equiv f$ and $(Q,Q_0)=(g,0) \equiv g$, with $f,g
\in \C$, equation (\ref{deriv}) gives
$$
[(d_*^{ \omega})^{(-E,0)}f, g]^{(0,1)} - [f,(d_*^{
\omega})^{(-E,0)}g]^{(0,1)}=(0,0),
$$
or, equivalently,
\begin{equation} \label{8.40}
[[(\Lambda,E),f]^{(0,1)}, g]^{(0,1)} -
[f,[(\Lambda,E),g]^{(0,1)}]^{(0,1)}=(0,0).
\end{equation}
The graded Jacobi identity for the bracket $[\cdot, \cdot]^{(0,1)}$
on $\Gamma(TM \times \R)$, ensures the validity of (\ref{8.40}).

Let us now take $(P,P_0)=(X,f) \in \Gamma(TM \times \R)$ and
$(Q,Q_0)=g \in \C$. Then,
\begin{eqnarray} \label{8.41}
\lefteqn{(d_*^{ \omega})^{(-E,0)}[(X,f),g]'^{(0,1)}= [ (\Lambda,E), [(X,f),g]^{(0,1)}]^{(0,1)} } \nonumber \\
& = & [(X,f),[ (\Lambda,E),g]^{(0,1)}]^{(0,1)}+[[(\Lambda,E),
(X,f)]^{(0,1)},g]^{(0,1)}\nonumber \\
& = & [(X,f),(d_*^{ \omega})^{(-E,0)}g]^{(0,1)} + [(d_*^{
\omega})^{(-E,0)}(X,f),g]^{(0,1)}- [((\Lambda,E)^\# \otimes
1)(d\omega,\omega)(X,f),g]^{(0,1)}\nonumber \\
& = & [(X,f),(d_*^{ \omega})^{(-E,0)}g]'^{(0,1)} + (\Lambda,E)^\#
\left( (d\omega,\omega)((X,f), (d_*^{ \omega})^{(-E,0)}g, \cdot)
\right)\nonumber \\
& & + [(d_*^{ \omega})^{(-E,0)}(X,f),g]'^{(0,1)}- [((\Lambda,E)^\#
\otimes
1)(d\omega,\omega)(X,f),g]^{(0,1)}\nonumber \\
& = & [(X,f),(d_*^{ \omega})^{(-E,0)}g]'^{(0,1)}+ [(d_*^{
\omega})^{(-E,0)}(X,f),g]'^{(0,1)},
\end{eqnarray}
which proves (\ref{deriv}) in this case.

When $(P,P_0)=(X,f)$ and $(Q,Q_0)=(Y,g)$ are two sections of $TM
\times \R$, equation (\ref{deriv}) is given by
\begin{equation} \label{8.42}
(d_*^{ \omega})^{(-E,0)}[(X,f),(Y,g)]'^{(0,1)}= [(d_*^{
\omega})^{(-E,0)}(X,f), (Y,g)]'^{(0,1)}+ [(X,f),(d_*^{
\omega})^{(-E,0)}(Y,g)]'^{(0,1)}. \\
\end{equation}
We compute,

\begin{eqnarray} \label{88.3}
\lefteqn{(d_*^{ \omega})^{(-E,0)}[(X,f),(Y,g)]'^{(0,1)}=
[(\Lambda,E),[(X,f),(Y,g)]^{(0,1)}]^{(0,1)} } \nonumber \\ &+ &
((\Lambda,E)^\# \otimes 1)(d\omega,\omega)([(X,f),(Y,g)]^{(0,1)})
\nonumber \\  & - &[(\Lambda,E),(\Lambda,E)^\#\left(
(d\omega,\omega)((X,f), (Y,g),
\cdot) \right)]^{(0,1)} \nonumber \\
& - &((\Lambda,E)^\# \otimes
1)(d\omega,\omega)((\Lambda,E)^\#((d\omega,\omega)((X,f), (Y,g),
\cdot))) \nonumber \\
& = & [(X,f), [(\Lambda,E), (Y,g)]^{(0,1)}]^{(0,1)} +
[[(\Lambda,E),(X,f)]^{(0,1)}, (Y,g)]^{(0,1)} \nonumber \\
&+ & ((\Lambda,E)^\# \otimes
1)(d\omega,\omega)([(X,f),(Y,g)]^{(0,1)}) \nonumber \\
& - & ((\Lambda,E)^\# \otimes
1)(d\omega,\omega)((\Lambda,E)^\#((d\omega,\omega)((X,f), (Y,g),
\cdot))) \nonumber \\
& + &(\Lambda,E)^\# \left( {\rm d}^{(0,1)}((d\omega,\omega)((X,f),
(Y,g), \cdot)) \right) +(\Lambda,E)^\#(d\omega,\omega) \left(
(d\omega,\omega)(X,f), (Y,g), \cdot) \right) \nonumber \\
\end{eqnarray}
where, in the last equality, we used (\ref{2.1}), the graded Jacobi
identity for the bracket $[\cdot, \cdot]^{(0,1)}$ and also the
following formula, that holds for any section $(\alpha,f)$ of $T^*M
\times \R$:
$$
[(\Lambda,E)^\#(\alpha ,f), (\Lambda,E)]= (\Lambda,E)^\#({\rm
d}^{(0,1)}(\alpha ,f)) + \frac{1}{2}
 [(\Lambda,E),(\Lambda,E)]^{(0,1)}((\alpha ,f)).
$$
On the other hand,

\begin{eqnarray} \label{88.4}
\lefteqn{[(d_*^{
\omega})^{(-E,0)}(X,f), (Y,g)]'^{(0,1)} + [(X,f),(d_*^{ \omega})^{(-E,0)}(Y,g)]'^{(0,1)}= } \nonumber \\
& = & [[(\Lambda,E),(X,f)]^{(0,1)},(Y,g)]^{(0,1)}\nonumber \\
& - & \left( \left((\Lambda,E)^\# \otimes
 ([(\Lambda,E),(X,f)]^{(0,1)})^\# +
 ([(\Lambda,E),(X,f)]^{(0,1)})^\# \otimes (\Lambda,E)^\# \right)
\otimes 1 \right)(d \omega, \omega)(Y,g)\nonumber \\
& + &[((\Lambda,E)^\# \otimes 1)(d \omega,
\omega) (X,f),(Y,g)]^{(0,1)} \nonumber \\
& - & \left( \left( (\Lambda,E)^\# \otimes \left( ((\Lambda,E)^\#
\otimes 1)(d \omega, \omega) (X,f) \right) ^\# \right. \right.
\nonumber
\\
&  & + \left. \left. \left(((\Lambda,E)^\# \otimes
1)(d \omega, \omega) (X,f) \right) ^\# \otimes (\Lambda,E)^\# \right) \otimes 1 \right)(d \omega, \omega)(Y,g)\nonumber \\
& + & [(X,f),[(\Lambda,E),(Y,g)]^{(0,1)}]^{(0,1)}\nonumber \\
& + & \left( \left((\Lambda,E)^\# \otimes
 ([(\Lambda,E),(Y,g)]^{(0,1)})^\# +
 ([(\Lambda,E),(Y,g)]^{(0,1)})^\# \otimes (\Lambda,E)^\# \right)
\otimes 1 \right)(d \omega, \omega)(X,f)\nonumber \\
& + & [(X,f),((\Lambda,E)^\# \otimes 1)(d \omega, \omega)
(Y,g)]^{(0,1)}\nonumber \\
& + & \left( \left( (\Lambda,E)^\# \otimes \left( ((\Lambda,E)^\#
\otimes 1)(d \omega, \omega) (Y,g) \right) ^\# \right. \right.
\nonumber
\\
&  & + \left. \left. \left(((\Lambda,E)^\# \otimes 1)(d \omega,
\omega) (Y,g) \right) ^\# \otimes (\Lambda,E)^\# \right) \otimes 1
\right)(d \omega, \omega)(X,f).
\end{eqnarray}
Comparing the terms of (\ref{88.3}) and (\ref{88.4}), we conclude,
after some computations, that (\ref{8.42}) holds if and only if, for
all $(\alpha, h), (\beta, l) \in \Gamma(T^*M \times \R)$,

\begin{eqnarray} \label{88.5}
\lefteqn{{\rm d}^{(0,1)}\left( (d \omega, \omega)((X,f), (Y,g),
\cdot)\right)((\Lambda,E)^\#(\alpha, h), (\Lambda,E)^\#(\beta, l) )}
\nonumber \\
& & - (d \omega, \omega)([(X,f),
(Y,g)]^{(0,1)},(\Lambda,E)^\#(\alpha, h), (\Lambda,E)^\#(\beta, l))\nonumber \\
& & + (d \omega, \omega)((Y,g),(\Lambda,E)^\#(\alpha, h),
([(\Lambda,E),(X,f)]^{(0,1)})^\#(\beta, l))\nonumber \\
& & -(d \omega, \omega)((Y,g),(\Lambda,E)^\#(\beta, l),
([(\Lambda,E),(X,f)]^{(0,1)})^\#(\alpha, h))\nonumber \\
& & + [(Y,g),((\Lambda,E)^\# \otimes 1)(d \omega,
\omega)(X,f)]^{(0,1)}((\alpha, h),(\beta, l))\nonumber \\
& & - (d \omega, \omega)((X,f),(\Lambda,E)^\#(\alpha, h),
([(\Lambda,E),(Y,g)]^{(0,1)})^\#(\beta, l))\nonumber \\
& & + (d \omega, \omega)((X,f),(\Lambda,E)^\#(\beta, l),
([(\Lambda,E),(Y,g)]^{(0,1)})^\#(\alpha, h))\nonumber \\
& & -[(X,f),((\Lambda,E)^\# \otimes 1)(d \omega,
\omega)(Y,g)]^{(0,1)}((\alpha, h),(\beta, l))=0.
\end{eqnarray}
After a long computation, we get that (\ref{88.5}) is equivalent to
$$\left( {\rm d}^{(0,1)} (d\omega,\omega) \right)((X,f), (Y,g),(\Lambda,E)^\#(\alpha, h), (\Lambda,E)^\#(\beta,
l))=0,$$

\noindent which holds since ${\rm d}^{(0,1)}
(d\omega,\omega)=(0,0)$.

\end{proof}

\section*{Acknowledgments} Research of J. M. Nunes da Costa supported by POCI/MAT/58452/2004 and CMUC-FCT.

\end{document}